\documentclass{daj}
\usepackage{amsmath}
\usepackage{amssymb}
 \numberwithin{equation}{subsection}
\newcommand{\per}{\mathrm{per\ }}
\newcommand{\haf}{\mathrm{haf\ }}
\newcommand{\PER}{\mathrm{PER\ }}
\newcommand{\UU}{\mathcal{U}}
\newtheorem{thm}{Theorem}[section]
\newtheorem{lem}{Lemma}[section]

\dajAUTHORdetails{%
  title = {Approximating Permanents and Hafnians}, 
  author = {Alexander Barvinok},
  plaintextauthor = {Alexander Barvinok},
  keywords = {permanent, hafnian, algorithm},
}   

\dajEDITORdetails{%
   year={2017},
   number={2},
   received={11 May 2016},   
   revised={20 December 2016},    
   published={13 January 2017},  
   doi={10.19086/da.1244},       
}   

\begin{document}

\begin{frontmatter}[classification=text]


\author[abar]{Alexander Barvinok\thanks{Supported by NSF Grant DMS 1361541.}}

\begin{abstract}
We prove that the logarithm of the permanent of an $n \times n$ real matrix $A$ and the logarithm of the hafnian of a $2n \times 2n$ real symmetric 
matrix $A$ can be approximated within an additive error $1 \geq \epsilon>0$ by a polynomial $p$ in the entries of $A$ of degree $O(\ln n -\ln \epsilon)$
provided the entries $a_{ij}$ of $A$ satisfy $\delta \leq a_{ij} \leq 1$ for an arbitrarily small $\delta >0$, fixed in advance.
Moreover, the polynomial $p$ can be computed in $n^{O(\ln n -\ln \epsilon)}$ time. 
 We also improve bounds for approximating $\ln \per A$, 
$\ln \haf A$ and logarithms of multi-dimensional permanents for complex matrices and tensors $A$.
\end{abstract}
\end{frontmatter}

\section{Main results: permanents}\label{sec1}

We discuss analytic methods of efficient approximation of permanents and hafnians of real and complex matrices as well as of their multi-dimensional versions,
objects of considerable interest in connection with problems in combinatorics \cite{LP09}, \cite{Mi78}, quantum physics \cite{AA13}, \cite{Ka16}, \cite{KK14} and computational complexity \cite{Va79}, \cite{J+04}. 
\bigskip

\subsection{Permanent}\label{subsec1.1} Let $A=\left(a_{ij}\right)$ be an $n \times n$ real or complex matrix. The {\it permanent} of $A$ is defined as 
$$\per A =\sum_{\sigma \in S_n} \prod_{i=1}^n a_{i \sigma(i)},
$$
where $S_n$ is the symmetric group of permutations of the set $\{1, \ldots, n\}$. It is a $\#P$-hard problem to compute the permanent of a given 0-1 matrix $A$ exactly \cite{Va79}, although a fully polynomial randomized approximation scheme is constructed for 
non-negative matrices \cite{J+04}. The permanent of an $n \times n$ non-negative matrix $A$ can be approximated within a factor of $e^n$ in deterministic polynomial time \cite{L+00} and the factor was improved to $2^n$ in \cite{GS14} (with a conjectured improvement to $2^{n/2}$).
 If one assumes that 
\begin{equation}\label{eq1.1.1}
\delta \ \leq \ a_{ij} \ \leq \ 1 \quad \textrm{for all} \quad i,j 
\end{equation}
and some $0 < \delta \leq 1$ fixed in advance, then the polynomial algorithm of \cite{L+00} actually results in an approximation factor of 
$n^{O(1)}$, where the implied constant in the ``$O$" notation depends on $\delta$, see also \cite{BS11}. Apart from that, deterministic polynomial time algorithms are known for special classes of matrices. For example, in \cite{GK10}, for any $\epsilon >0$, fixed in advance, a polynomial time algorithm is constructed to approximate $\per A$ within a factor of $(1+\epsilon)^n$ if $A$ is the adjacency matrix of a constant degree expander. 
We also note that in \cite{F+04} a simple randomized algorithm is shown to approximate $\per A$ within a subexponential in $n$ factor provided 
(\ref{eq1.1.1}) holds with some $0 < \delta \leq 1$, fixed in advance.

In this paper, we present a quasi-polynomial deterministic algorithm, which, given an $n \times n$ matrix $A=\left(a_{ij}\right)$ satisfying (\ref{eq1.1.1}) with 
some $0 < \delta \leq 1$, fixed in advance, and an $\epsilon >0$ approximates $\per A$ within a relative error $\epsilon$ in 
$n^{O(\ln n - \ln \epsilon)}$ time. The implied constant in the ``$O$" notation depends on $\delta$. 

More precisely, we prove the following result.

\begin{thm}\label{th1.2} For any $0 < \delta \leq 1$ there exists $\gamma=\gamma(\delta) >0$ such that for any positive integer $n$ and 
any $0< \epsilon < 1$ there exists a polynomial $p=p_{n, \delta, \epsilon}$ in the entries $a_{ij}$ of an $n \times n$ matrix $A$ such that $\deg p \leq \gamma\left(\ln n - \ln \epsilon\right)$ and
$$\left| \ln \per A - p(A) \right|\leq \epsilon$$
for all $n \times n$ real matrices $A=\left(a_{ij}\right)$ satisfying 
$$\delta \ \leq \ a_{ij} \ \leq \ 1 \quad \textrm{for all} \quad i, j.$$
\end{thm}
We show that the polynomial $p$ can be computed in quasi-polynomial time $n^{O(\ln n -\ln \epsilon)}$, where the implied constant in the ``$O$" notation depends on $\delta$ alone.

Our approach continues a line of work started in \cite{Ba16} and continued in \cite{Ba15}, \cite{BS16} and \cite{Re15}.
The main idea is to relate approximability of a polynomial with its {\it complex} zeros. For a complex number $z=a+ib$, we denote 
by $\Re\thinspace z =a$ and $\Im\thinspace z=b$, the real and imaginary parts of $z$ correspondingly. In what follows, we always choose the standard branch of 
$\arcsin x$, $\arccos x$ and $\arctan x$ for real $x$, so that 

\begin{equation*}
\begin{split} -\frac{\pi}{2}  \ \leq \ &\arcsin x \ \leq \ \frac{\pi}{2} \quad \text{for} \quad  -1 \ \leq \ x \ \leq \ 1,\\
0 \ \leq \ &\arccos x \ \leq \ \pi \quad \text{for} \quad -1 \ \leq \ x \ \leq \ 1 \qquad \text{and} \\
-\frac{\pi}{2} \  < \ &\arctan x \ < \ \frac{\pi}{2} \quad \text{for} \quad x \in {\mathbb R}.
\end{split}
\end{equation*}

We deduce Theorem \ref{th1.2} from the following result. 

\begin{thm}\label{th1.3}
Let us fix a real $0 \leq  \eta < 1$ and let
$$\tau=(1-\eta) \sin\left(\frac{\pi}{4} - \arctan \eta \right)  \ > \ 0.$$
Let $Z=\left(z_{ij}\right)$ be an $n \times n$ complex matrix such that 
$$\left|1 - \Re\thinspace z_{ij} \right| \ \leq \ \eta \quad \text{and} \quad \left| \Im\thinspace z_{ij} \right| \ \leq \ \tau \quad \text{for all} 
\quad 1 \leq i, j \leq n.$$
Then $\per Z \ne 0$.
\end{thm}

There is an interest in computing permanents of {\it complex} matrices \cite{AA13}, \cite{KK14}, \cite{Ka16}.  Ryser's algorithm, see for example Chapter 7 of \cite{Mi78}, computes the permanent of an $n \times n$ matrix $A$ over an arbitrary field exactly in $O(n 2^n)$ time. Exact polynomial time algorithms are known for 
rather restricted classes of matrices, such as matrices of a fixed rank \cite{Ba96} or matrices for which the support of non-zero entries is a graph of a fixed tree-width \cite{CP16}.
In \cite{Fu00}, a randomized polynomial time algorithm is constructed which computes the permanent of a complex matrix within a (properly defined) relative error $\epsilon >0$ in 
$O\left(3^{n/2} \epsilon^{-2}\right)$ time. In \cite{Gu05}, a randomized algorithm is constructed which approximates $\per A$ for a complex $n \times n$ 
matrix $A$ within an additive error $\epsilon \|A\|^n$, where $\|A\|$ is the operator norm of $A$, in time polynomial in $n$ and $1/\epsilon$, see also \cite{AA13} for an exposition.

In this paper, we prove the following results.

\begin{thm}\label{th1.4}
Let $Z=\left(z_{ij}\right)$ be an $n \times n$ complex matrix such that 
$$\left| 1- z_{ij} \right| \ \leq \ 0.5 \quad \text{for all} \quad 1 \leq i, j \leq n.$$
Then $\per Z \ne 0$.
\end{thm}

Since $\per Z \ne 0$, we can choose a branch of $\ln \per Z$ when the conditions of Theorem \ref{th1.4} are satisfied (for convenience, we always choose 
the branch for which $\ln \per Z$ is real if $Z$ is a real matrix).
We deduce from Theorem \ref{th1.4} the following approximation result.

\begin{thm}\label{th1.5}
For every $0 \leq \eta < 0.5$ there exists a constant $\gamma=\gamma(\eta)>0$ such that  for every positive integer $n$ and every real $0 < \epsilon < 1$ there exists a polynomial $p=p_{n, \eta, \epsilon}$ in the entries of an $n \times n$ complex matrix $A=\left(a_{ij}\right)$ such that 
$\deg p \leq \gamma\left(\ln n - \ln \epsilon\right)$ and 
$$\left| \ln \per A - p(A) \right| \ \leq \ \epsilon$$
for $n \times n$ complex matrices $A=\left(a_{ij}\right)$ satisfying 
$$\left| 1- a_{ij} \right| \ \leq \ \eta \quad \text{for all} \quad i, j.$$ 
\end{thm}
Moreover, the polynomial $p$ can be computed in $n^{O(\ln n -\ln \epsilon)}$ time, where the implied constant in the ``$O$" notation depends on $\eta$ alone.

A version of Theorem \ref{th1.4} with a weaker bound of $0.195$ instead of $0.5$ and a more complicated proof was obtained in \cite{Ba16}.
Theorem \ref{th1.5} is also implicit in \cite{Ba16}. We present its proof here since it serves as a stepping stone for the proof of Theorem \ref{th1.2}.

It is not known whether the bound $0.5$ in Theorems \ref{th1.4} and \ref{th1.5} can be increased, although one can show (see Section \ref{sec4}) that it cannot be increased
to $\sqrt{2}/2 \approx 0.707$.

Theorems \ref{th1.4} and \ref{th1.5} state, roughly, that the permanent behaves nicely as long as the matrix is not too far in the $\ell^{\infty}$-distance from from the matrix $J$ of all 1s. Applied to an arbitrary $n \times n$ positive matrix $A$, Theorem \ref{th1.5} implies that $\per A$ can be approximated deterministically within a relative error $0 < \epsilon < 1$
in quasi-polynomial time $n^{O(\ln n -\ln \epsilon)}$ as long as the entries of $A$ are within some multiplicative factor $\gamma < 3$, fixed in advance, of 
each other. 

Let $A$ be an $n \times n$ complex matrix such that the $\ell^{\infty}$-distance from $A$ to 
the complex hypersurface of $n \times n$ matrices $Z$ satisfying $\per Z=0$ is at least $\delta_0>0$. It follows from our proof that for any $0 < \delta < \delta_0$ and any 
$0 < \epsilon < 1$ there exists a polynomial $p_A$ in the entries of an $n\times n$ matrix such that $\left| \ln \per B - p_A(B) \right| \leq \epsilon$ 
for any matrix $B$ within distance $\delta$ in the $\ell^{\infty}$-distance from $A$ and $\deg p_A =O(\ln n -\ln \epsilon)$, where the implied constant in the 
``$O$" notation depends only on $\delta$ and $\delta_0$. However, for a general $A \ne J$, finding the polynomial $p_A$ may be computationally hard.

 Theorems \ref{th1.2} and \ref{th1.3} are of a different nature: there we allow the entries $a_{ij}$ to be arbitrarily close to $0$ but insist that the imaginary part of $a_{ij}$ get smaller as $a_{ij}$ approach $0$. Theorem \ref{th1.2} implies that for a positive $n \times n$ matrix $A$, the value $\per A$ can be approximated 
 deterministically within a relative error $0 < \epsilon < 1$ in quasi-polynomial time $n^{O(\ln n -\ln \epsilon)}$ as long as the entries of $A$ are within some multiplicative factor 
 $\gamma \geq 1$, arbitrarily large, but fixed in advance, of each other. It follows from our proofs that a similar to Theorem \ref{th1.2} approximation result holds for complex matrices $A=\left(a_{ij}\right)$ with $\delta \leq \Re\thinspace a_{ij} \leq 1$ and $\left| \Im a_{ij}\right| \leq \tau_0$
 for some fixed $\tau_0=\tau_0(\delta) >0$.
 
 So far, we approximated permanents of real or complex matrices that are close to the matrix $J$ of all 1s in the $\ell^{\infty}$-distance. Next, we consider 
 matrices that are close to $J$ in the maximum $\ell^1$-distance over all rows and columns. 

 \begin{thm}\label{th1.6}
 Let $\alpha \approx 0.278$ be the real solution of the equation $\alpha e^{1+\alpha}=1$.
 Let $Z=\left(z_{ij}\right)$ be an $n \times n$ complex matrix such that 
 $$\sum_{j=1}^n \left| 1- z_{ij}\right| \ \leq \ \frac{\alpha n}{4}  \quad \text{for} \quad i=1, \ldots, n$$ 
and 
$$\sum_{i=1}^n \left| 1 - z_{ij} \right| \ \leq \ \frac{\alpha n}{4}  \quad \text{for} \quad j=1, \ldots, n.$$
Then $\per Z \ne 0$. 
 \end{thm}
 
 Since $\per Z \ne 0$, we can choose a branch of $\ln \per Z$ when the conditions of Theorem \ref{th1.6} are satisfied.
 We obtain the following result. 
 
 \begin{thm}\label{th1.7}
  For every $0 \leq \eta < \alpha/4$, where $\alpha \approx 0.278$ is the constant in Theorem \ref{th1.6}, there exists a constant 
 $\gamma=\gamma(\eta) >0$ such that for every positive integer $n$ and every real $0 < \epsilon < 1$ there exists a polynomial $p=p_{n, \eta, \epsilon}$ in the entries 
 of an $n \times n$ matrix $A=\left(a_{ij}\right)$ such that $\deg p \leq \gamma(\ln n -\ln \epsilon)$ and 
 $$\left| \ln \per A -p(A) \right| \ \leq \ \epsilon$$
 for $n \times n$ complex matrices $A=\left(a_{ij}\right)$ satisfying 
 $$\sum_{j=1}^n \left| 1- a_{ij}\right| \ \leq \ \eta n \quad \text{for} \quad i=1, \ldots, n$$ 
 and 
 $$\sum_{i=1}^n \left| 1- a_{ij}\right| \ \leq \ \eta n \quad \text{for} \quad j=1, \ldots, n.$$
 \end{thm}
 
 Again, the polynomial $p_{n, \eta, \epsilon}$ can be constructed in $n^{O(\ln n -\ln \epsilon)}$ time, where the implied constant in the ``$O$" notation depends on 
 $\eta$ only. Note that Theorem \ref{th1.7} is applicable to 0-1 matrices having not too many (not more than $7\%$) zeros in every row and column as well as to real matrices with some positive and some negative entries. It is not known whether the bound in Theorems \ref{th1.6} and \ref{th1.7} are optimal.

 \section{Main results: hafnians}\label{sec2}
 
Some of our results immediately extend from permanents to hafnians.

\subsection{Hafnian}\label{subsec2.1} Let $A=\left(a_{ij}\right)$ be a $2n \times 2n$ symmetric real or complex matrix. The {\it hafnian} of $A$ is defined as
$$\haf A=\sum_{\{i_1, j_1\}, \ldots, \{i_n, j_n\}} a_{i_1 j_1} \cdots a_{i_n j_n},$$
where the sum is taken over $(2n)!/2^n n!$ unordered partitions of the set $\{1, \ldots, 2n\}$ into $n$ pairwise disjoint unordered pairs 
$\{i_1, j_1\}, \ldots, \{i_n, j_n\}$, see for example, Section 8.2 of \cite{Mi78}. Just as the permanent of the biadjacency matrix of a bipartite graph enumerates the perfect matchings in the graph, the hafnian of the adjacency matrix of a graph enumerates the perfect matchings in the graph.
In fact, for any $n \times n$ matrix $A$ we have 
$$\haf \left( \begin{matrix} 0 & A \\ A^T & 0 \end{matrix} \right) =\per A,$$
and hence computing the permanent of an $n \times n$ matrix reduces to computing the hafnian of a symmetric $2n \times 2n$ matrix. 

Computationally, the hafnian appears to be a more complicated object than the permanent. No fully polynomial (randomized or deterministic) approximation scheme is known to compute the hafnian of a non-negative symmetric matrix and no deterministic polynomial time algorithm to approximate the hafnian of a $2n \times 2n$ non-negative matrix within an exponential factor of $c^n$ for some absolute constant $c>1$ is known (though there is a randomized polynomial time algorithm achieving such an approximation \cite{Ba99}, see also \cite{R+16} for cases when the algorithm approximates within a subexponential factor). On the other hand, if the entries $a_{ij}$ of the matrix $A=\left(a_{ij}\right)$ satisfy (\ref{eq1.1.1}) for some $\delta >0$, fixed in advance, there is a polynomial time algorithm approximating $\haf A$ within a factor of $n^{O(1)}$, where the implied constant in the ``$O$" notation depends on $\delta$ \cite{BS11}. 

In this paper, we prove the following versions of Theorem \ref{th1.2} and \ref{th1.3}.

\begin{thm}\label{th2.2}
 For any $0 < \delta \leq 1$ there exists $\gamma=\gamma(\delta) >0$ such that for any positive integer $n$ and 
any $0< \epsilon < 1$ there exists a polynomial $p=p_{n, \delta, \epsilon}$ in the entries $a_{ij}$ of a $2n \times 2n$ symmetric matrix $A$ such that 
$\deg p \leq \gamma(\ln n -\ln \epsilon)$ and 
$$\left| \ln \haf A - p(A) \right|\leq \epsilon$$
for all $2n \times 2n$ real symmetric matrices $A=\left(a_{ij}\right)$ satisfying 
$$\delta \ \leq \ a_{ij} \ \leq \ 1 \quad \text{for all} \quad i, j.$$
\end{thm}
The polynomial $p_{n, \delta, \epsilon}$ can be computed in $n^{O(\ln n -\ln \epsilon)}$ time, where the implied constant in the ``$O$" notation depends on 
$\delta$ alone. Consequently, we obtain a deterministic quasi-polynomial algorithm to approximate the hafnian of a positive matrix $A=\left(a_{ij}\right)$ satisfying (\ref{eq1.1.1}) within any given relative error $\epsilon >0$.

As is the case with permanents, we deduce Theorem \ref{th2.2} from a result on the complex zeros of the hafnian.
\begin{thm}\label{th2.3}
Let us fix a real $0 \leq  \eta <1$ and and let 
$$\tau \ = \ (1-\eta) \sin\left(\frac{\pi}{4} -\arctan \eta \right) >0.$$ 
Let $Z=\left(z_{ij}\right)$ be an $2n \times 2n$ symmetric complex matrix such that 
$$\left|1 -\Re\thinspace z_{ij}\right| \ \leq \ \eta \quad \text{and} \quad \left| \Im\thinspace z_{ij} \right| \ \leq \ \tau \quad \text{for all} 
\quad 1 \leq i, j \leq n.$$
Then $\haf Z \ne 0$.
\end{thm}

We also obtain the following versions of Theorems \ref{th1.4} and \ref{th1.5}.

\begin{thm}\label{th2.4}
 Let $Z=\left(z_{ij}\right)$ be an $2n \times 2n$ symmetric complex matrix such that 
$$\left| 1- z_{ij} \right| \ \leq \ 0.5 \quad \text{for all} \quad i, j.$$
Then $\haf Z \ne 0$. 
\end{thm}

As before, for matrices $Z$ satisfying the condition of Theorem \ref{th2.4}, we choose a branch of $\ln \haf Z$ in such a way so that $\ln \haf Z$ is real if 
$Z$ is a real matrix. We obtain the following result.

\begin{thm}\label{th2.5}
For any $0 \leq \eta < 0.5$ there exists $\gamma=\gamma(\eta) >0$ and for any positive integer $n$ and real $0 < \epsilon < 1$
there exists a polynomial $p=p_{n, \eta, \epsilon}$ in the entries of $2n \times 2n$ complex symmetric matrix $A=\left(a_{ij}\right)$ such that 
$\deg p \leq \gamma (\ln n -\ln \epsilon)$ and 
$$\left| \ln \haf A - p(A) \right| \ \leq \ \epsilon$$
provided 
$$|1 -a_{ij}| \ \leq \ \eta \quad \text{for all} \quad i, j.$$
\end{thm}

As before, the polynomial $p_{n, \eta, \epsilon}$ can be computed in $n^{O(\ln n -\ln \epsilon)}$ time, where the implied constant in the ``$O$" notation 
depends on $\eta$ alone.

 Our approach can be extended to a variety of 
{\it partition functions} \cite{Ba15}, \cite{BS16}. In Section \ref{sec3}, we show how to extend it to multi-dimensional permanents of 
tensors.

\section{Main results: multi-dimensional permanents}\label{sec3}

\subsection{Multi-dimensional permanent}\label{subsec3.1}
Let $A=\left(a_{i_1 \ldots i_d}\right)$ be a $d$-dimensional $n \times \ldots \times n$ array (tensor) filled with 
$n^d$ real or complex matrices. We define the {\it permanent} of $A$ by 
$$\PER A =\sum_{\sigma_2, \ldots, \sigma_d \in S_n} \prod_{i=1}^n a_{i \sigma_2(i) \ldots \sigma_d(i)}.$$

In particular, if $d=2$ then $A$ is an $n \times n$ matrix and $\PER A=\per A$. If $d \geq 3$ it is an NP-complete problem 
to tell $\PER A$ from $0$ if $A$ is a tensor with 0-1 entries, since the problem reduces to finding whether a given $d$-partite hypergraph has 
a perfect matching. 

We define a {\it slice} of $A$ as the array of $n^{d-1}$ entries of $A$ with one of the indices $i_1, \ldots, i_d$ fixed to a particular value and the remaining $(d-1)$ indices varying arbitrarily. Hence $A$ has altogether $nd$ slices. If $d=2$ and $A$ is a matrix then a slice is a row or a column.

We note that for $d >2$ there are several different notions of the permanent of a tensor, cf., for example, \cite{LL14}.

We obtain the following extension of Theorem \ref{th1.4}.
\begin{thm}\label{th3.2}
For an integer $d \geq 2$, let us choose 
$$\eta_d= \sin \frac{\theta}{2} \cos \frac{(d-1) \theta}{2}$$
for some $\theta =\theta_d >0$ such that $(d-1) \theta \ < \ 2\pi/3$.
Hence $0 < \eta_d < 1$ and we can choose $\eta_2 =0.5$, $\eta_3 = \sqrt{6}/9 \approx 0.272$, $\eta_4 \approx 0.184$
and $\eta_d=\Omega\left(\frac{1}{d}\right)$.

Let $Z=\left(z_{i_1 \ldots i_d}\right)$ be a $d$-dimensional complex $n \times \ldots \times n$ array such that 
$$\left|1-z_{i_1 \ldots i_d} \right| \ \leq \ \eta_d \quad \text{for all} \quad i_1 \ldots i_d.$$
Then $\PER Z \ne 0$.
\end{thm}

A version of Theorem \ref{th3.2} with weaker bounds $\eta_2=0.195$, $\eta_3=0.125$ and $\eta_4=0.093$ and a more complicated proof was obtained in
\cite{Ba16}.

Since $\PER Z \ne 0$, we can choose a branch of $\ln \PER Z$ when the conditions of Theorem \ref{th3.2} are satisfied (as before, we choose the branch for which $\ln \PER Z$ is real if $Z$ is a real tensor).
As a corollary, we obtain the following approximation result.

\begin{thm}\label{th3.3}
 For an integer $d \geq 2$, let us choose $0 \leq \eta < \eta_d$, where $\eta_d$ is the constant in Theorem \ref{th3.2}. Then there exists
$\gamma=\gamma(d, \eta) >0$ and for every integer $n$ and real $0 < \epsilon < 1$ there exists a polynomial $p=p_{d, \eta, \epsilon, n}$ in the entries of a $d$-dimensional $n \times \ldots \times n$ complex tensor $A=\left(a_{i_1 \ldots i_d}\right)$ such that 
$\deg p \leq \gamma(\ln n -\ln \epsilon)$ and 
$$\left| \ln \PER A - p(A) \right| \ \leq \ \epsilon$$
provided 
$$\left|1 -a_{i_1 \ldots i_d} \right| \ \leq \ \eta \quad \text{for all} \quad 1 \leq i_1, \ldots, i_d \leq n.$$
\end{thm}

The polynomial $p_{d, \eta, \epsilon, n}$ can be computed in $n^{O(\ln n -\ln \epsilon)}$ time, where the implied constant in the ``$O$" notation depends only on 
$d$ and $\eta$.

While we were unable to obtain exact equivalents of Theorems \ref{th1.2} and \ref{th2.2}, our approach produces the following approximation result for multi-dimensional permanents. 

\begin{thm}\label{th3.4}
For an integer $d \geq 2$, let 
$$\eta_d = \tan \frac{\pi}{4(d-1)}$$
so that $\eta_2=1$,  $\eta_3 = \sqrt{2}-1 \approx 0.414$, $\eta_4= 2-\sqrt{3} \approx 0.268$, etc.

For any $0 \leq \eta < \eta_d$ there is a constant $\gamma=\gamma(d, \eta)$ and for any positive integer $n$ and real $0 < \epsilon < 1$ there is a polynomial 
$p=p_{d, \eta, \epsilon, n}$ is the entries of a $d$-dimensional $n \times \cdots \times n$ tensor such that $\deg p \leq \gamma(\ln n -\ln \epsilon)$ and 
$$\left| \ln \PER A - p(A) \right| \ \leq \ \epsilon$$ 
for any $d$-dimensional $n \times \ldots \times n$ real tensor $A=\left(a_{i_1 \ldots i_d}\right)$ satisfying 
$$\left| 1 - a_{i_1 \ldots i_d} \right| \ \leq \ \eta \quad \text{for all} \quad 1 \leq i_1, \ldots, i_d \leq n.$$
\end{thm}

Again, the polynomial $p_{d, \eta, \epsilon, n}$ can be computed in $n^{O(\ln n -\ln \epsilon)}$ time, where the implied constant in the ``$O$" notation depends only on $d$ and $\eta$. For example, for $n \times n \times n$ tensors $A$ with positive real entries, we obtain a quasi-polynomial algorithm to approximate $\PER A$ if the entries of $A$ are within a factor $\gamma < \sqrt{2}+1 \approx 2.414$ of each other. Note that Theorem \ref{th3.3} for $n \times n \times n$ tensors $A$ with positive real entries guarantees the existences of a quasi-polynomial algorithm to approximate $\PER A$ if the entries of $A$ are within a factor of $(1+\sqrt{6}/9)/(1-\sqrt{6}/9) \approx 1.748$ of each other.

As before, the proof is based on the absence of zeros of $\PER A$ in a particular domain. Namely, we deduce Theorem \ref{th3.4} from the following result.

\begin{thm}\label{th3.5}
For an integer $d \geq 2$, let 
$\eta_d$ be the constant of Theorem \ref{th3.4}.
Let us fix a real $0 \leq \eta < \eta_d$ and let 
$$\tau=(1-\eta) \sin \left( \frac{\pi}{4(d-1)} -\arctan  \eta \right) > 0.$$
Let $Z=\left(z_{i_1 \ldots i_d}\right)$  be a $d$-dimensional tensor of complex numbers such that 
$$\left|1- \Re\thinspace z_{i_1 \ldots i_d} \right| \ \leq \ \eta \quad \text{and} \quad \left| \Im\thinspace z_{i_1 \ldots i_d} \right| \ \leq \ \tau$$
for all $1 \leq i_1, \ldots, i_d \leq n$.

Then $\PER Z \ne 0$.
\end{thm}

Finally, we obtain multi-dimensional versions of Theorems \ref{th1.6} and \ref{th1.7}.  

\begin{thm}\label{th3.6}
 Let $\alpha \approx 0.278$ be the real solution of the equation 
$\alpha e^{1+\alpha}=1$. For an integer $d \geq 2$, let 
$$\eta_d=\frac{\alpha^{d-1} (d-1)^{d-1}}{d^d}.$$
 Let $Z=\left(z_{i_1 \ldots i_d}\right)$ 
be a $d$-dimensional complex $n \times \ldots \times n$ array such that the sum of 
$\left|1-z_{i_1 \ldots i_d}\right|$ over each slice of $Z$ does not exceed $\eta_d n^{d-1}$.

Then $\PER Z \ne 0$.
\end{thm}

In other words, $\PER Z \ne 0$ if each slice of $Z$ is sufficiently close to the 
array of 1s in the $\ell^1$-distance. While for each fixed $d$, the allowed distance is of the order of $n^{d-1}$, it decreases exponentially with $d$, unlike the allowed $\ell^{\infty}$-distance in Theorem \ref{th3.2}, which decreases as $1/d$.

We obtain the following corollary.

\begin{thm}\label{th3.7}
 For every integer $d \geq 2$ and every
$0 \ \leq \ \eta \ < \ \eta_d$, 
where $\eta_d$ is the constant of Theorem \ref{th3.6}, there exists a constant $\gamma=\gamma(d, \eta) >0$ such that for any positive integer $n$ and real $0 < \epsilon < 1$ there is a polynomial $p=p_{d,\eta,\epsilon,n}$ in the entries of a $d$-dimensional $n \times \ldots \times n$ tensor such that $\deg p \leq \gamma(\ln n -\ln \epsilon)$ and
$$\left| \ln \PER A - p(A) \right| \ \leq \ \epsilon$$
for any $d$-dimensional $n \times \ldots \times n$ tensor $A=\left(a_{i_1 \ldots i_d}\right)$ for which the sum of $\left|1-a_{i_1 \ldots i_d}\right|$ over each slice of $A$ does not exceed $\eta n^{d-1}$.
\end{thm}

Again, the polynomial $p_{d,\eta,\epsilon,n}$ can be computed in $n^{O(\ln n -\ln \epsilon)}$ time, where the implied constant in the ``$O$" notation depends on 
$d$ and $\eta$ alone. Theorem \ref{th3.7} is applicable to 0-1 tensors $A$, which contain a small (and exponentially decreasing with $d$) fraction of 0s in each slice.
\bigskip

In Section \ref{sec4}, we prove Theorems \ref{th1.4}, \ref{th2.4} and \ref{th3.2}.

In Section \ref{sec5}, we prove Theorems \ref{th1.3}, \ref{th2.3} and \ref{th3.4}.

In Section \ref{sec6}, we prove Theorems \ref{th1.6} and \ref{th3.6}.

In Section \ref{sec7}, we prove Theorems \ref{th1.5}, \ref{th1.7}, \ref{th2.5}, \ref{th3.3} and \ref{th3.7}.

In Section \ref{sec8}, we prove Theorems \ref{th1.2}, \ref{th2.2} and \ref{th3.4}.

Finally, in Section \ref{sec9}, we discuss possible ramifications and open questions.

\section{Proofs of Theorems \ref{th1.4}, \ref{th2.4} and \ref{th3.2}}\label{sec4}

We start with a simple geometric argument regarding angles between non-zero complex numbers. We identify ${\mathbb C}={\mathbb R}^2$, thus identifying complex numbers with vectors in the plane. We denote by $\langle \cdot, \cdot \rangle$ the standard inner product in ${\mathbb R}^2$, so that 
$$\langle a, b \rangle = \Re \left(a \overline{b}\right) \quad \text{for} \quad a, b \in {\mathbb C}$$ and by 
$| \cdot |$ the corresponding Euclidean norm (the modulus of a complex number).

\begin{lem}\label{le4.1}
 Let $0 \leq \theta < 2\pi/3$ be real and let $u_1, \ldots, u_n \in {\mathbb C}$ be non-zero 
complex numbers such that the angle between any two $u_i$ and $u_j$ does not exceed $\theta$. Let 
$$u=u_1+ \ldots + u_n.$$
Then 
\begin{enumerate}
\item We have 
$$|u| \ \geq \ \left(\cos \frac{\theta}{2}\right) \sum_{j=1}^n \left| u_j \right|.$$
\item Let $\alpha_1, \ldots, \alpha_n$ and $\beta_1, \ldots, \beta_n$ be complex numbers 
such that 
$$\left| 1- \alpha_j \right| \ \leq \ \eta \quad \text{and} \quad \left| 1-\beta_j \right| \ \leq \ \eta$$
for some 
$$0 \ \leq \ \eta \ < \ \cos \frac{\theta}{2}$$
and $j=1, \ldots, n$.
Let 
$$v = \sum_{j=1}^n \alpha_j u_j \quad \text{and} \quad w =\sum_{j=1}^n \beta_j u_j.$$
Then $v \ne 0$, $w \ne 0$ and the angle between $v$ and $w$ does not exceed 
$$ 2 \arcsin \frac{\eta}{\cos(\theta/2)}.$$
\end{enumerate}
\end{lem}
\textbf{Proof.} Part (1) and its proof is due to Boris Bukh \cite{Bu15}. If $0$ is in the convex hull of $u_1, \ldots, u_n$ then, by the 
Carath\'eodory Theorem, we conclude that $0$ is in the convex hull of some three vectors $u_i, u_j$ and $u_k$ and hence the angle between some two vectors $u_i$ and $u_j$ is at least $2\pi/3$, which is a contradiction. Therefore, $0$ is not in the convex hull of $u_1, \ldots, u_n$ and hence the vectors $u_1, \ldots, u_n$ lie in a cone $K \subset {\mathbb C}$ of measure at most $\theta$ with vertex at $0$. 

Let us consider the orthogonal projection of each vector $u_j$ onto the bisector of $K$. Then the length of the projection of $u_j$ 
is at least $|u_j| \cos (\theta/2)$ and hence the length of the orthogonal projection of $u$ onto the bisector of $K$ is at least 
$$\left(\cos \frac{\theta}{2}\right) \sum_{j=1}^n |u_j|.$$
Since the length of $u$ is at least as large as the length of its orthogonal projection, the proof of Part (1) follows.

To prove Part (2), we note that 
$$|v-u| = \left| \sum_{j=1}^n \left(\alpha_j -1 \right) u_j \right| \ \leq \ \eta \sum_{j=1}^n |u_j|.$$
From Part (1), we conclude that $|v-u| < |u|$. Therefore, $v=(v-u)+u \ne 0$ and the angle between $v$ and $u$ does not exceed 
$$\arcsin \frac{|v-u|}{|u|} \ \leq \ \arcsin \frac{\eta}{\cos (\theta/2)}.$$
Similarly, $w=(w-u)+u \ne 0$ and the angle between $w$ and $u$ does not exceed 
$$\arcsin \frac{|w-u|}{|u|} \ \leq \ \arcsin \frac{\eta}{\cos (\theta/2)}.$$
Therefore, the angle between $v$ and $w$ does not exceed 
$$2 \arcsin \frac{\eta}{\cos (\theta/2)}$$
and the proof of Part (2) follows.
$\hfill \square$

\subsection{Proof of Theorem \ref{th1.4}}\label{subsec4.2} For a positive integer $n$, let $\UU_n$ be the set of $n \times n$ complex 
matrices $Z=\left(z_{ij}\right)$ such that 
\begin{equation}\label{eq4.2.1}
\left| 1- z_{ij} \right| \ \leq \ 0.5 \quad \textrm{for all} \quad i, j.
\end{equation}
We prove by induction on $n$ the following statement: 
\bigskip

For any $Z \in \UU_n$ we have $\per Z \ne 0$ and, moreover, if $A, B \in \UU_n$ are two matrices that differ in one row (or in one column) only then the angle between non-zero complex numbers $\per A$ and $\per B$ does not exceed $\pi/2$.
\bigskip

The statement obviously holds for $n=1$. Assuming that the statement holds for matrices in $\UU_{n-1}$ with $n \geq 2$, let us consider two matrices $A, B \in \UU_n$ that differ in one row or in one column only. Since the permanent of a matrix does not change when the rows or columns of the matrix are permuted or when the matrix is transposed, without loss of generality we assume that $B$ is obtained from $A$ by replacing the entries 
$a_{1j}$ of the first row by complex numbers $b_{1j}$ for $j=1, \ldots, n$. Let $A_j$ be the $(n-1) \times (n-1)$ matrix obtained from $A$ by crossing out the first row and the $j$-th column. Then 
\begin{equation}\label{eq4.2.2}
\per A = \sum_{j=1}^n a_{1j} \per A_j \quad \text{and} \quad \per B=\sum_{j=1}^n b_{1j} \per A_j.
\end{equation}
We observe that $A_j \in \UU_{n-1}$ for $j=1, \ldots, n$ and, moreover, any two matrices $A_{j_1}$ and $A_{j_2}$ after a suitable permutation of columns differ in one column only. Hence by the induction hypothesis, we have $\per A_j \ne 0$ for $j=1, \ldots, n$ and the angle between any two non-zero complex numbers $\per A_{j_1}$ and $\per A_{j_2}$ does not exceed $\pi/2$. 
Applying Part (2) of Lemma \ref{le4.1}
with 
$$\theta =\frac{\pi}{2},\quad \eta=\frac{1}{2}, \quad u_j =\per A_j, \quad
  \alpha_j = a_{1j} \quad \text{and} \quad  \beta_j =b_{1j} \quad \text{for} \quad j=1, \ldots, n,$$
 we conclude that 
$\per A \ne 0$, $\per B \ne 0$ and the angle between non-zero complex numbers $\per A$ and $\per B$ does not exceed
$$2 \arcsin \frac{0.5}{\cos(\pi/4)} = 2 \arcsin \frac{\sqrt{2}}{2} = \frac{\pi}{2},$$
which concludes the induction step.
\hfill $\square$

One can observe that $\eta=0.5$ is the largest value of $\eta$ for which the equation 
$$\theta =2 \arcsin \frac{\eta}{\cos(\theta/2)}$$
has a solution $\theta < 2\pi/3$ and hence the induction in Section \ref{subsec4.2}  can proceed. It is not known whether the constant $0.5$  in Theorem \ref{th1.4} can be increased.
Since 
$$\per A=0 \quad \text{where} \quad A= \left( \begin{matrix} \frac{1+i}{2} & \frac{1-i}{2} \\ \frac{1-i}{2} & \frac{1+i}{2} \end{matrix} \right)=0,$$
the value of $0.5$ in Theorem \ref{th1.4} cannot be replaced by $\sqrt{2}/2 \approx 0.707$. Moreover, as Boris Bukh noticed \cite{Bu15}, 
we have 
$$\per \left( A \otimes J_m\right) =0,$$
where $A$ is a matrix as above, $m$ is odd and $J_m$ is an $m \times m$ matrix filled with 1s. 

\subsection{Proof of Theorem \ref{th2.4}}\label{subsec4.3}
 The proof is very similar to that of Section \ref{subsec4.2}. For a positive integer $n$, we define $\UU_{n}$ as the set of $2n \times 2n$ symmetric complex matrices $Z=\left(z_{ij}\right)$ satisfying (\ref{eq4.2.1})
  and prove by induction on $n$ that for any $Z \in \UU_n$ we have $\haf Z \ne 0$ and if $A, B \in \UU_n$ are two matrices that differ only in the $k$-th row and in the $k$-th column for some unique $k$ then the angle between non-zero complex numbers $\haf A$ and $\haf B$ does not exceed $\pi/2$. 
  
The statement obviously holds for $n=1$. Suppose that $n >1$. Since the hafnian of the matrix does not change under a simultaneous permutation of rows and columns, without loss of generality we may assume that $A$ and $B$ differ in the first row and first column only.
Instead of the Laplace expansion (\ref{eq4.2.2}), we use the recurrence
\begin{equation}\label{eq4.3.1}
\haf A=\sum_{j=2}^{2n} a_{1j} \haf A_j \quad \text{and} \quad \haf B =\sum_{j=2}^{2n} b_{1j} \haf A_j 
\end{equation}
where $A_j$ is the $(2n-2) \times (2n-2)$ matrix obtained from $A$ by crossing out the first row and the first column and the $j$-th row and the $j$-th column. We observe that, up to a simultaneous permutation of rows and columns, any two matrices $A_{j_1}$ and $A_{j_2}$ differ only in the $k$-th row and $k$-th column for some $k$ and the induction proceeds as in Section \ref{subsec4.2}.
\hfill $\square$

\subsection{Proof of Theorem \ref{th3.2}}\label{subsec4.4} By and large, the proof proceeds as in Section \ref{subsec4.2}. For a positive integer $n$, we define $\UU_n$ as the set of $n \times \ldots \times n$ complex arrays 
$Z=\left(z_{i_1 \ldots i_d}\right)$ such that 
$$\left|1 - z_{i_1 \ldots i_d} \right| \ \leq \ \eta_d \quad \text{for all} \quad i_1, \ldots, i_d.$$
We prove by induction on $n$ the following statement:

\bigskip
For any $Z \in \UU_n$ we have $\PER Z \ne 0$ and, moreover, if $A, B \in \UU_n$ are two tensors that differ in one slice only, then the angle between non-zero complex numbers $\PER A$ and $\PER B$ does not exceed $\theta$.
\bigskip 

If $n=1$ then the angle between $\PER A$ and $\PER B$ does not exceed
$$2 \arcsin \eta_d \ < \ 2 \arcsin \left(\sin \frac{\theta}{2}\right) = \theta$$ and 
the statement holds. Assuming that $n \geq 2$, let us consider two tensors $A, B \in \UU_n$ that differ in one slice only.  Without loss of generality, we assume that $B$ is obtained from $A$ by replacing the ``top slice" numbers 
$a_{1i_2 \ldots i_d}$ with numbers $b_{1 i_2 \ldots i_d}$. We use a $d$-dimensional version of the Laplace expansion:
\begin{equation}\label{eq4.4.1}
\begin{split} \PER A =&\sum_{1 \leq i_2, \ldots, i_d \leq n} a_{1i_2 \ldots i_d} \PER A_{i_2 \ldots i_d}  \quad \text{and} \\
\PER B=&\sum_{1 \leq i_2, \ldots i_d \leq n} b_{1 i_2 \ldots i_d} \PER A_{i_2 \ldots i_d},
\end{split}
\end{equation}
where $A_{i_2 \ldots i_d}$ is the $(n -1) \times \ldots \times (n-1)$ tensor obtained from $A$ by crossing out the $d$ slices obtained by fixing the first index to 1, the second index to $i_2$, $\ldots$, the last index to $i_d$. It remains to notice that any two 
tensors $A_{i_2 \ldots i_d}$ and $A_{i_2' \ldots i_d'}$ differ in at most $d-1$ slices, and hence by the induction hypothesis 
we have $\PER A_{i_2 \ldots i_d} \ne 0$, $\PER A_{i_2' \ldots i_d'} \ne 0$ and the angle between the two non-zero complex numbers  does not exceed $(d-1) \theta$. Applying Part (2) of Lemma \ref{le4.1}, we conclude that $\PER A \ne 0$, $\PER B \ne 0$ and the angle between non-zero complex numbers $\PER A$ and $\PER B$ does not exceed 
$$2 \arcsin\frac{\eta_d}{\cos\left( \frac{(d-1) \theta}{2}\right)} =\theta,$$
which completes the induction.
\hfill $\square$

\section{Proofs of Theorems \ref{th1.3}, \ref{th2.3} and \ref{th3.5}}\label{sec5}

As in Section \ref{sec4}, we start with a simple geometric lemma.
\begin{lem}\label{le5.1}
 Let $u_1, \ldots, u_n \in {\mathbb C}$ be non-zero complex numbers such that the angle between any two $u_i$ and $u_j$ does not exceed $\pi/2$. Let
$$v=\sum_{j=1}^n \alpha_j u_j \quad \text{and} \quad w=\sum_{j=1}^n \beta_j u_j$$
for some complex numbers $\alpha_1, \ldots, \alpha_n$ and $\beta_1, \ldots, \beta_n$. 
\begin{enumerate}
\item Suppose that $\alpha_1, \ldots, \alpha_n$ are non-negative real and that $\beta_1, \ldots, \beta_n$ are real such that 
$$\left| \beta_j \right| \ \leq \ \alpha_j \quad \text{for} \quad j=1, \ldots, n$$
Then $|w| \ \leq \ |v|$.
\item Suppose that $\alpha_1, \ldots, \alpha_n$ and $\beta_1, \ldots, \beta_n$ are real such that 
$$\left| 1 - \alpha_j \right| \ \leq \ \eta \quad \text{and} \quad \left| 1- \beta_j \right| \ \leq \ \eta \quad \text{for} \quad j=1, \ldots, n$$
for some $0 \leq \eta < 1$.
Then $v \ne 0$, $w \ne 0$ and the angle between $v$ and $w$ does not exceed $2 \arctan \eta$.
\item Suppose that 
\begin{equation*}
\begin{split} &\left|1 -\Re\thinspace \alpha_j \right| \ \leq \ \eta, \quad \left|1 - \Re\thinspace \beta_j \right| \ \leq \ \eta \quad \text{and} \\
&\left| \Im\thinspace \alpha_j \right| \ \leq \ \tau, \quad \left| \Im\thinspace \beta_j \right| \ \leq \ \tau \quad \text{for} \quad j=1, \ldots, n 
\end{split}
\end{equation*}
for some $0 \leq \eta < 1$ and some $0 \leq \tau < 1-\eta$. Then $v \ne 0$, $w \ne 0$ and the angle between $v$ and $w$ does not exceed
$$2 \arctan \eta + 2 \arcsin \frac{\tau}{1-\eta}.$$
\end{enumerate}
\end{lem}
\textbf{Proof.} Since 
$$\langle u_i, u_j \rangle \ \geq \ 0 \quad \text{for all} \quad i, j,$$
in Part (1) we obtain 
$$|w|^2 =\sum_{1 \leq i, j\leq n} \beta_i \beta_j \langle u_i, u_j \rangle \ \leq \ \sum_{1 \leq i, j \leq n} \alpha_i \alpha_j \langle u_i, u_j \rangle = |v|^2$$
and the proof of Part (1) follows.

To prove Part (2), let 
$$u=\frac{v + w}{2} =\sum_{j=1}^n \left(\frac{\alpha_j + \beta_j}{2}\right) u_j \quad \text{and} \quad x=\frac{v -w}{2} =
\sum_{j=1}^n \left(\frac{\alpha_j -\beta_j}{2}\right) u_j,$$
so that
 $$v=u+x \quad \text{and} \quad  w=u-x.$$ For $j=1, \ldots, n$, we have 
$$\eta \left(\alpha_j + \beta_j\right) - \left(\alpha_j - \beta_j \right) = \beta_j (1+\eta) - \alpha_j (1-\eta) \ \geq \ (1-\eta)(1+\eta) - (1+\eta)(1-\eta) 
\ \geq \ 0,$$ 
from which it follows that 
$$\left| \frac{\alpha_j - \beta_j}{2} \right| \ \leq \ \eta \left(\frac{\alpha_j + \beta_j}{2}\right) \quad \text{for} \quad j=1, \ldots, n$$ 
and hence by Part (1) we have 
$$|x| \leq \eta |u|.$$ It follows that $v \ne 0$, $w \ne 0$ and that the angle between $v$ and $w$ is 
$$\arccos \frac{\langle v, w \rangle}{|v| |w|}.$$
We have
$$\langle v, w \rangle = \langle u+x, u-x \rangle = |u|^2 - |x|^2 \ > \ 0$$
and 
$$|v|^2 + |w|^2 = \langle u+x, u +x \rangle + \langle u-x, u-x \rangle = 2 |u|^2 + 2|x|^2,$$
so that 
$$|v||w| \ \leq \ |u|^2 + |x|^2$$
with the equality attained when $|v|^2=|w|^2 = |u|^2 + |x|^2$ and $x$ is orthogonal to $u$. Hence for given $|u|$ and $|x|$ the largest angle of 
$$ \arccos \frac{|u|^2 - |x|^2}{|u|^2 + |x|^2}$$
 between $v$ and $w$ is attained when $x$ is orthogonal to $u$ and is equal to 
 $$2 \arctan\frac{|x|}{|u|} \ \leq \ 2 \arctan \eta,$$
 which completes the proof of Part (2).

 To prove Part (3), let 
 \begin{equation*}
 \begin{split} v'=&\sum_{j=1}^n \left(\Re\thinspace \alpha_j\right) u_j, \quad v''=\sum_{j=1}^n \left(\Im\thinspace \alpha_j\right) u_j, \quad w'=\sum_{j=1}^n \left(\Re\thinspace \beta_j\right) u_j \\ \text{and} \quad  w''=&\sum_{j=1}^n \left(\Im\thinspace \beta_j\right) u_j. 
 \end{split}
 \end{equation*}
 By Part (2), $v' \ne 0$, $w' \ne 0$ and the angle between $v'$ and $w'$ does not exceed $\theta=2\arctan \eta$. 
 Since
 $$\Re\thinspace \alpha_j, \ \Re\thinspace \beta_j \ \geq \ 1-\eta \quad \textrm{and} \quad \left| \Im\thinspace \beta_j \right|, \ \left| \Im\thinspace \beta_j \right| 
 \ \leq \ \tau
 \quad \text{for} \quad j=1, \ldots, n,$$
 from Part (1), we conclude that
 $$\left|v''\right| \ \leq \ \frac{\tau}{1-\eta} \left| v' \right| \quad \text{and} \quad \left| w''\right| \ \leq \ \frac{\tau}{1-\eta}\left| w'\right|.$$
 Since $0 \leq \tau < 1-\eta$, we have $v=v'+i v'' \ne 0$, $w=w'+i w'' \ne 0$ and the angle between $v$ and $v'$ and the angle between $w$ and $w'$ do not exceed 
 $$\omega=\arcsin \frac{ \tau}{1-\eta}.$$
 Therefore the angle between $v$ and $w$ does not exceed $\theta + 2\omega$ and the proof of Part (3) follows. 
\hfill $\square$

\subsection{Proof of Theorem \ref{th1.3}}\label{subsec5.2}
For a positive integer $n$, let $\UU_n$ be the set of $n \times n$ complex matrices $Z=\left(z_{ij}\right)$ such that 
\begin{equation}\label{eq5.2.1}
\left| 1- \Re\thinspace z_{ij} \right|  \ \leq \ \eta  \quad \text{and} \quad \left| \Im\thinspace z_{ij} \right| \ \leq \ \tau \quad \text{for all} \quad i, j. 
\end{equation}
We prove by induction on $n$ the following statement:
\bigskip

For any $Z \in \UU_n$ we have $\per Z \ne 0$ and, moreover, if $A, B \in \UU_n$ are two matrices that differ in one row (or in one column) only, then 
the angle between non-zero complex numbers $\per A$ and $\per B$ does not exceed $\pi/2$.
\bigskip

Since $\tau <  1-\eta$, the statement holds for $n=1$. Assuming that the statement holds for matrices in $\UU_{n-1}$ with $n \geq 2$, let us consider two matrices $A, B \in \UU_n$ 
that differ in one row or in one column only. As in Section \ref{subsec4.2}, without loss of generality we assume that $B$ is obtained from $A$ by replacing the entries 
$a_{1j}$ of the first row by complex numbers $b_{1j}$ for $j=1, \ldots, n$. Let $A_j$ be the $(n-1) \times (n-1)$ matrix obtained from $A$ by crossing out the first row and the $j$-th column. 
We observe that $A_j \in \UU_{n-1}$ for $j=1, \ldots, n$ and, moreover, any two matrices $A_{j_1}$ and $A_{j_2}$ after a suitable permutation of columns differ in one column only. Hence by the induction hypothesis, we have $\per A_j \ne 0$ for $j=1, \ldots, n$ and the angle between any two non-zero complex numbers $\per A_{j_1}$ and $\per A_{j_2}$ does not exceed $\pi/2$.  Using the Laplace expansion (\ref{eq4.2.2}) and applying Part (3) of Lemma \ref{le5.1} with 
$$u_j =\per A_j,\quad \alpha_j=a_{1j} \quad \text{and} \quad  \beta_j=b_{1j} \quad \text{for} \quad j=1, \ldots, n,$$
 we conclude that $\per A \ne 0$, $\per B \ne 0$ and that the angle between non-zero complex numbers $\per A$ and $\per B$ does not exceed 
 $$2 \arctan \eta + 2 \arcsin \frac{\tau}{1-\eta} =\frac{\pi}{2},$$
  which completes the induction.
\hfill $\square$

\subsection{Proof of Theorem \ref{th2.3}}\label{subsec5.3} The proof is very similar to that of Section \ref{subsec5.2}. For a positive integer $n$, we define $\UU_{n}$ as the set of $2n \times 2n$ symmetric complex matrices $Z=\left(z_{ij}\right)$ satisfying (\ref{eq5.2.1}) and prove by induction on $n$ that for any $Z \in \UU_n$ we have $\haf Z \ne 0$ and if $A, B \in \UU_n$ are two matrices that differ only in the $k$-th row and in the $k$-th column for some unique $k$ then the angle between non-zero complex numbers $\haf A$ and $\haf B$ does not exceed $\pi/2$. 

Since $\tau < 1-\eta$, the statement holds for $n=1$. Suppose that $n > 1$. As in Section \ref{subsec4.3}, without loss of generality we assume that $A$ and $B$ differ in the first row and column only.  Let $A_j$ be the $(2n-2) \times (2n-2)$ matrix obtained from $A$ by crossing out the first row and the first column and the $j$-th row and the $j$-th column. As in Section \ref{subsec4.3},  we observe that, up to a simultaneous permutation of rows and columns (which does not change the hafnian), any two matrices $A_{j_1}$ and $A_{j_2}$ differ only in the $k$-th row and $k$-th column for some $k$. Using the expansion (\ref{eq4.3.1}), we complete the induction as in Section \ref{subsec5.2}.
\hfill $\square$

\subsection{Proof of Theorem \ref{th3.5}}\label{subsec5.4} By and large, the proof proceeds as in Section \ref{subsec5.2}.  For a positive integer $n$, we define $\UU_n$ as the set of $n \times \ldots \times n$ complex arrays 
$Z=\left(z_{i_1 \ldots i_d}\right)$ such that 
$$\left| 1- \Re\thinspace z_{i_1 \ldots i_d} \right| \ \leq \ \eta \quad \text{and} \quad \left|\Im\thinspace z_{i_1 \ldots i_d}\right| \ \leq \ \tau$$
for all $1 \leq i_1, \ldots, i_d \leq n$. We prove by induction on $n$ the following statement:
\bigskip

For any $Z \in \UU_n$ we have $\PER Z \ne 0$ and, moreover, if $A, B \in \UU_n$ are two tensors that differ in one slice only, then the angle between non-zero complex numbers $\PER A$ and $\PER B$ does not exceed $\frac{\pi}{2(d-1)}$.
\bigskip

Since $\tau < 1-\eta$, the statement holds for $n=1$. Assuming that $n \geq 2$, let us consider two tensors $A, B \in \UU_n$ that differ in one slice only. 
As in Section \ref{subsec4.4}, we assume that $B$ is obtained from $A$ by replacing the top slice numbers $a_{1i_2 \ldots i_d}$ with numbers $b_{1 i_2 \ldots i_d}$ and define 
the $(n -1) \times \ldots \times (n-1)$ tensor $A_{i_2 \ldots i_d}$ as the tensor obtained from $A$ by crossing out the $d$ slices obtained by fixing the first index to 1, the second index to $i_2$, $\ldots$, the last index to $i_d$. As in Section \ref{subsec4.4}, any two tensors $A_{i_2 \ldots i_d}$ and $A_{i_2' \ldots i_d'}$ differ in at most $d-1$ slices, and hence by the induction hypothesis 
we have $\PER A_{i_2 \ldots i_d} \ne 0$, $\PER A_{i_2' \ldots i_d'} \ne 0$ and the angle between the two non-zero complex numbers  does not exceed
$\pi/2$. Using the $d$-dimensional version (\ref{eq4.4.1}) of the Laplace expansion and Part (3) of Lemma \ref{le5.1}, we conclude that $\PER A \ne 0$, $\PER B \ne 0$ and the angle between non-zero complex numbers $\PER A$ and 
$\PER B$ does not exceed 
$$2 \arctan \eta + 2 \arcsin \frac{\tau}{1-\eta}=\frac{\pi}{2(d-1)},$$
which completes the induction.
\hfill $\square$

\section{Proofs of Theorems \ref{th1.6} and \ref{th3.6}}\label{sec6}

Since Theorem \ref{th1.6} is a particular case of Theorem \ref{th3.6} for $d=2$, we prove the latter theorem. We use a combinatorial interpretation of the multi-dimensional permanent in terms of matchings in a hypergraph.

\subsection{The matching polynomial of a hypergraph}\label{subsec6.1} Let us fix an integer $d \geq 2$. Let $V$ be a finite set and let 
$E \subset \binom{V}{d}$ be 
a family of $d$-subsets of $V$. The pair $H=(V, E)$ is called a $d$-hypergraph with set $V$ of vertices and set $E$ of edges. An unordered set $e_1, \ldots, e_k$ of pairwise disjoint edges of $H$ is called a {\it matching} (we agree that the empty set of edges is a matching). Given a map $w: E \longrightarrow {\mathbb C}$ that assigns complex weights $w(e)$ to the edges $e \in E$ of $H$, we define the {\it weight} of a matching $e_1, \ldots, e_k$ as the the product $w(e_1) \cdots w(e_k)$ of weights of the edges of the matching. We agree that the weight of the empty matching is 1. We define the {\it matching polynomial} as the sum of weights of all 
matchings (including the empty one) in $H$:
$$P_H(w) = \sum_{\substack{e_1, \ldots, e_k \\ \textrm{is a matching}}}  w(e_1) \cdots w(e_k).$$

\begin{lem}\label{le6.2}
 Let $H=(V, E)$ be a $d$-hypergraph and let $w: E \longrightarrow {\mathbb C}$ be complex weights on its edges. 
Suppose that 
$$\sum_{\substack{ e \in E:\\ v \in e}} |w(e)| \ \leq \ \frac{(d-1)^{d-1}}{d^d}  \quad \textrm{for all} \quad v \in V.$$
Then $P_H(w) \ne 0$.
\end{lem}
\textbf{Proof.} For a set $S \subset V$ of vertices, we denote by $H-S$ the hypergraph with set $V \setminus S$ of vertices and set $E' \subset E$ of edges that do not contain vertices from $S$. Abusing notation, we denote the restriction of weights $w: E \longrightarrow {\mathbb C}$ onto $E'$ also by $w$.
 We prove by induction on the number $|V|$ of vertices that $P_H(w) \ne 0$ and, moreover, for every vertex $v \in V$ we have 
\begin{equation}\label{eq6.2.1}
\left| 1 - \frac{P_{H-\{v\}}(w)}{P_H(w)} \right| \ \leq \ \frac{1}{d-1}. 
\end{equation}
If $|V| < d$ then $H$ has no edges and hence $P_H(w)=P_{H-v}(w)=1$. Suppose now that $|V| \geq d$. We observe the following recurrence:
\begin{equation}\label{eq6.2.2}
P_H(w)=P_{H-\{v\}}(w) + \sum_{\substack{ e \in E:\\ v \in e}} w(e) P_{H-e}(w), 
\end{equation}
where the $P_{H-\{v\}}(w)$ accounts for the matchings in $H$ not containing $v$ and the sum accounts for the matchings of $H$ containing $v$. By the induction 
hypothesis, $P_{H-\{v\}}(w) \ne 0$, so we rewrite (\ref{eq6.2.2}) as 
\begin{equation}\label{eq6.2.3}
\frac{P_H(w)}{ P_{H-\{v\}}(w)} = 1+ \sum_{\substack{e \in E:\\ v \in e}} w(e) \frac{P_{H-e}(w)}{P_{H-\{v\}}(w)}. 
\end{equation}
If there are no edges $e$ containing $v$ then $P_H(w)=P_{H-\{v\}}(w)$ and (\ref{eq6.2.1}) follows. Otherwise, let $e=\{v, v_2, \ldots, v_d\}$ be an edge containing 
$v$. Telescoping, we obtain
\begin{equation}\label{eq6.2.4}
\frac{P_{H-e}(w)}{P_{H-\{v\}}(w)}=\frac{P_{H-e}(w)}{P_{H-\{v, v_2 \ldots, v_{d-1}\}}(w) } \frac{P_{H-\{v, v_2, \ldots, v_{d-1}\}}(w)}{P_{H-\{v, v_2 \ldots, v_{d-2}\}}(w)} 
\cdots \frac{P_{H-\{v, v_2\}}(w)}{P_{H-\{v\}}(w)}. 
\end{equation}
By the induction hypothesis, each ratio in the right hand side of (\ref{eq6.2.4}) does not exceed $d/(d-1)$ in the absolute value, and hence 
$$\left| \frac{P_{H-e}(w)}{P_{H-\{v\}}(w)}\right| \ \leq \ \left(\frac{d}{d-1}\right)^{d-1}.$$
Therefore, from (\ref{eq6.2.3}) we obtain 
\begin{equation}\label{eq6.2.5}
\left| 1- \frac{P_H(w)}{P_{H-\{v\}}(w)}\right| \ \leq \ \frac{(d-1)^{d-1}}{d^d} \left(\frac{d}{d-1}\right)^{d-1} = \frac{1}{d}, 
\end{equation}
from which it follows that $P_H(w) \ne 0$. Denoting 
$$z = \frac{P_{H-\{v\}}(w)}{P_{H}(w)},$$ from (\ref{eq6.2.5}) we have a chain of implications
\begin{equation*}
\begin{split}
 \left| 1- \frac{1}{z} \right| \ \leq \ \frac{1}{d} \quad \Longrightarrow \quad  &\left| \frac{z-1}{z} \right| \ \leq \ \frac{1}{d} \quad \Longrightarrow \quad \left| \frac{z}{z-1} \right| \ \geq \ d \\ \Longrightarrow \quad &\left| \frac{1}{z-1} \right| \ \geq \ d-1 \quad \Longrightarrow \quad
|1-z| \ \leq \ \frac{1}{d-1} 
\end{split}
\end{equation*}
proving (\ref{eq6.2.1}).
\hfill $\square$

The bound of Lemma \ref{le6.2} and to some extent its proof agrees with those of \cite{HL72} for the roots of the matching polynomial of a graph.

Next, we need a weaker version on an estimate from \cite{Wa03}. 
\begin{lem}\label{le6.3}
Let $\alpha \approx 0.278$ be the constant of Theorem \ref{th3.6}, so that $\alpha e^{1+\alpha} =1$. For a positive integer $n$, let 
$$p_n(z)=\sum_{k=0}^n \frac{z^k}{k!}.$$
Then 
$$p_n(z) \ne 0 \quad \text{provided} \quad |z| \ \leq \ \alpha n.$$
\end{lem}
\textbf{Proof.} We observe that if $|z| \leq \alpha$ then 
$$\left|z e^{1-z} \right| \ \leq \ |z| e^{1+|z|} \ \leq \ 1$$
and hence
\begin{equation*}
\begin{split}
\left| 1- e^{-nz} p_n(nz)\right|= &\left| e^{-nz} \sum_{k=n+1}^{\infty} \frac{(nz)^k}{k!} \right| =
\left| \left( z e^{1-z}\right)^n e^{-n} \sum_{k=n+1}^{\infty} \frac{n^k z^{k-n}}{k!} \right| \\
\ \leq \ &e^{-n} \sum_{k=n+1}^{\infty} \frac{n^k}{k!} \ < \ 1, 
\end{split}
\end{equation*}
so that $p_n(nz) \ne 0$.
\hfill $\square$

Finally, we need a theorem of Szeg\H{o}, see for example, Chapter IV of \cite{Ma66} and also \cite{BB09} for generalizations.
\begin{thm}\label{th6.4}
 Let 
$$f(z)=\sum_{k=0}^n a_k z^k \quad \text{and} \quad g(z)=\sum_{k=0}^n b_k z^k$$
be complex polynomials. We define the Schur product $h=f \ast g$ by 
$$h(z) = \sum_{k=0}^n c_k z^k \quad \text{where} \quad c_k=\frac{a_k b_k}{\binom{n}{k}} \quad \text{for} \quad k=0, \ldots, n.$$
Suppose that $f(z) \ne 0$ whenever $|z| \leq r_1$ and $g(z) \ne 0$ whenever $|z| \leq r_2$ for some $r_1, r_2 >0$. 

Then 
$h(z) \ne 0$ whenever $|z| \leq r_1 r_2$.
\end{thm}

\subsection{Proof of Theorem \ref{th3.6}}\label{subsec6.5} Let $H$ be the complete $d$-partite hypergraph with set $V$ of $nd$ vertices, split into $d$ parts and vertices in each part numbered $1$ through $n$. Each edge of $H$ consist of exactly one vertex from each part and we let the weight of edge $\left(i_1, \ldots, i_d\right)$ equal to 
$w_{i_1 \ldots i_d}=z_{i_1 \ldots i_d}-1$. For $k=1, \ldots, n$, let  $W_k$ be the total weight of all matchings in $H$ consisting of exactly $k$ edges. We write
\begin{equation*}
\begin{split}
 \PER Z = &\sum_{\sigma_2, \ldots, \sigma_d \in S_n} \prod_{i=1}^n z_{i \sigma_2(i) \ldots \sigma_d(i)} =
\sum_{\sigma_2, \ldots, \sigma_d \in S_n}  \prod_{i=1}^n \left(1 + w_{i \sigma_2(i) \ldots \sigma_d(i)} \right) \\
=&\sum_{\sigma_2, \ldots, \sigma_d \in S_n}\left(1+ \sum_{k=1}^n \sum_{1 \leq i_1 < \ldots <  i_k \leq n}  w_{i_1 \sigma_2(i_1) \ldots \sigma_d(i_1)}
\cdots w_{i_k \sigma_2(i_k) \ldots \sigma_d(i_k)}\right) \\
=&\sum_{k=0}^n \left((n-k)!\right)^{d-1} W_k. 
\end{split}
\end{equation*}
Let us define a univariate polynomial 
$$f(z)=\sum_{k=0}^n W_k z^k.$$
Then $f(z)$ is the value of the matching polynomial $P_H$ on the scaled weights $z w_{i_1 \ldots i_d}$ and from Lemma \ref{le6.2} we conclude that 
\begin{equation}\label{eq6.5.1}
f(z) \ne 0 \quad \text{provided} \quad |z| \ \leq \ \frac{1}{(\alpha n)^{d-1}}. 
\end{equation}
Let $p_n$ be the polynomial of Lemma \ref{le6.3}. Applying Lemma \ref{le6.3} and Theorem \ref{th6.4} to the Schur product $h=f \ast p_n \ast \cdots \ast p_n$ of 
$f$ and $d-1$ polynomials $p_n$, we conclude from (\ref{eq6.5.1}) that
$$ h(z)=\left(\frac{1}{n!}\right)^{d-1} \sum_{k=0}^n ((n-k)!)^{d-1} W_k z^k  \ne 0 \quad \text{provided} \quad |z| \leq 1.$$
In particular, $h(1) \ne 0$ and hence $\PER Z \ne 0$.
\hfill $\square$

\section{Proofs of Theorems \ref{th1.5}, \ref{th1.7}, \ref{th2.5}, \ref{th3.3} and \ref{th3.7}}\label{sec7}

We need the following simple result first obtained in \cite{Ba16}. For completeness, 
we give its proof here. 

\begin{lem}\label{le7.1}
Let $g: {\mathbb C} \longrightarrow {\mathbb C}$ be a polynomial and let $\beta >1$ be real such that $g(z) \ne 0$ for all $|z| \leq \beta$. Let us choose a branch of 
$$f(z) =\ln g(z) \quad \text{for} \quad |z| \leq \beta$$
and let 
$$T_m(z)=f(0) + \sum_{k=1}^m  \frac{f^{(k)}(0)}{k!} z^k$$
be the Taylor polynomial of $f(z)$ of degree $m$ computed at $z=0$. Then 
$$\left| f(1) - T_m(1) \right| \ \leq \ \frac{\deg g}{(m+1) \beta^m (\beta-1)}.$$
\end{lem}
\textbf{Proof.}
Without loss of generality, we assume that $n=\deg g > 0$. Let $z_1, \ldots, z _n \in {\mathbb C}$ be the roots of $g$, each root is listed with its multiplicity. 
Hence we can write 
$$g(z) =g(0) \prod_{j=1}^n \left(1 - \frac{z}{z_j}\right) \quad \text{where} \quad |z_j| > \beta \quad \text{for} \quad j=1, \ldots, n$$
and 
$$f(z) =f(0) + \sum_{j=1}^n \ln \left( 1- \frac{z}{z_j} \right) \quad \text{for all} \quad |z| \leq 1.$$
Using the Taylor series expansion for the logarithm, we obtain 
$$\ln \left(1-\frac{1}{z_j}\right) =-\sum_{k=1}^m \frac{1}{k z_j^k} + \xi_j $$
where 
$$\left| \xi_j \right| = \left| - \sum_{k=m+1}^{\infty} \frac{1}{k z_j^k}  \right| \ \leq \ \frac{1}{m+1} \sum_{k=m+1}^{\infty} \frac{1}{ \beta^k}  \ = \ \frac{1}{(m+1)\beta^m (\beta-1)}.$$
Since
$$T_m(1) =f(0)- \sum_{j=1}^n \sum_{k=1}^m \frac{1}{k z_j^k},$$
the proof follows.
\hfill $\square$

It follows from Lemma \ref{le7.1} that as long as the roots of a polynomial $g(z)$ stay at distance at least $\beta$ away from $0$ for some fixed  $\beta >1$, then to approximate $\ln g(1)$ within an additive error $\epsilon$, we can use the Taylor polynomial of $f(z)=\ln g(z)$ at $z=0$ of degree $m=O(\ln \deg g - \ln \epsilon)$, where the implied constant in the ``$O$" notation depends on $\beta$ only.

\subsection{Computing the derivatives}\label{subsec7.2} As is discussed in \cite{Ba16}, the computation of the first $m$ derivatives
$f^{(1)}(0), \ldots, f^{(m)}(0)$ of $f(z)=\ln g(z)$ reduces to the computation of the first $m$ derivatives 
$g^{(1)}(0), \ldots, g^{(m)}(0)$ of $g$. Indeed, 
$$f^{(1)}(z) =\frac{g^{(1)}(z)}{g(z)} \quad \text{and hence} \quad g^{(1)}(z)= f^{(1)}(z) g(z).$$
Therefore,
\begin{equation*}
\begin{split}
 g^{(2)}(z) =&f^{(2)}(z) g(z) + f^{(1)}(z) g^{(1)}(z), \\ g^{(3)}(z) = &f^{(3)}(z) g(z) + 2 f^{(2)}(z) g^{(1)}(z) + f^{(1)}(z) g^{(2)}(z)
\end{split}
\end{equation*}
and
\begin{equation}\label{eq7.2.1}
g^{(k)}(0) =\sum_{j=0}^{k-1} \binom{k-1}{j} g^{(j)}(0) f^{(k-j)}(0)
\end{equation}
where $g^{(0)}(0)=g(0) \ne 0$. Writing equations (\ref{eq7.2.1}) for $k=1, \ldots, m$ we obtain a non-singular triangular system of linear equations in $f^{(k)}(0)$ with numbers 
$g(0)\ne 0$ on the diagonal 
from which the values of $f^{(1)}(0), \ldots, f^{(m)}(0)$ can be computed in $O(m^2)$ time from the values of $g(0), g^{(1)}(0), \ldots, g^{(m)}(0)$.
Thus
\begin{equation*}
\begin{split}
 f^{(1)}(0)=&\frac{g^{(1)}(0)}{g(0)}, \quad f^{(2)}(0)=\frac{g^{(2)}(0)-f^{(1)}(0) g^{(1)}(0)}{g(0)}=\frac{g^{(2)}(0)}{g(0)} -\frac{\left(g^{(1)}(0)\right)^2}{(g(0))^2} \\
f^{(3)}(0)= &\frac{g^{(3)}(0)-2f^{(2)}(0)g^{(1)}(0)-f^{(1)}(0) g^{(2)}(0)}{g(0)}\\=&\frac{g^{(3)}(0)}{g(0)} -\frac{3g^{(2)}(0) g^{(1)}(0)}{(g(0))^2} 
+\frac{2\left(g^{(1)}(0)\right)^3}{(g(0))^3},
\end{split}
\end{equation*}
and, generally, $f^{(k)}(0)$ is a linear combination of expressions of the 
type
$$\frac{g^{(k_1)}(0) \cdots g^{(k_s)}(0)}{(g(0))^p} \quad \text{where} \quad k_1 + \ldots + k_s =k \quad \text{and} \quad p \geq 1$$
with integer coefficients. 

Note that computing $f^{(k)}(0)$ from $g^{(k)}(0)$ is akin to computing cumulants of a distribution from its moments.

\subsection{Proof of Theorem \ref{th1.5}}\label{subsec7.3}
Let $J=J_n$ be the $n \times n$ matrix filled with 1s and let $A=\left(a_{ij}\right)$ be an $n \times n$ complex matrix satisfying the conditions of the theorem. We define a univariate polynomial 
$$g(z) = \per \bigl(J + z (A-J)\bigr) \quad \text{for} \quad z \in {\mathbb C}.$$
so that $\deg g \leq n$,  
$$g(0) = \per J =n! \quad \text{and} \quad g(1)=\per A.$$ 
Moreover, by Theorem \ref{th1.4} we have 
$$g(z) \ne 0 \quad \text{provided} \quad |z| \leq \beta \quad \text{where} \quad \beta = \frac{0.5}{\eta} > 1.$$
Let us choose the branch of 
$$f(z) =\ln g(z) \quad \text{for} \quad |z| \leq \beta$$
so that $f(0)$ is real
and let 
$$T_m(z) =f(0) + \sum_{k=1}^m \frac{f^{(k)}(0)}{k!} z^k$$
be the Taylor polynomial of degree $m$ computed at $z=0$. It follows from Lemma \ref{le7.1} that for some constant $\gamma=\gamma(\eta)>0$ and 
 integer $m \leq \gamma(\ln n -\ln \epsilon)$ we have 
$$\left| \ln \per A - T_m(1)\right| = \left| f(1) - T_m(1) \right| \ \leq \ \epsilon.$$

It remains to show that $T_m(1)$ is a polynomial $p$ in the entries $a_{ij}$ of the matrix $A$ of degree at most $m$. In view of Section \ref{subsec7.2} and the fact that 
$g(0)=n!$, it suffices to check that $g^{(k)}(0)$ is a polynomial in the entries $a_{ij}$ of the matrix $A$ of degree at most $k$ which can be computed 
in $n^{O(k)}$ time, where the implied constant in the ``$O$" notation is absolute.
We have 
$$g(z) =\sum_{\sigma \in S_n} \prod_{i=1}^n \left(1+ z\left(a_{i \sigma(i)} -1 \right)\right)$$
and hence for $k \geq 1$
$$g^{(k)}(0)=\sum_{\sigma \in S_n} \sum_{(i_1, \ldots, i_k)} \left(a_{i_1\sigma(i_1)}-1 \right) \cdots \left(a_{i_k \sigma(i_k)}-1\right),$$
where the last sum is taken over all ordered sets $(i_1, \ldots, i_k)$ of distinct numbers between $1$ and $n$. By symmetry, we can further write
$$g^{(k)}(0)=(n-k)! \sum_{\substack{(i_1, \ldots, i_k) \\ (j_1, \ldots, j_k)}} \left(a_{i_1j_1} -1\right) \cdots \left(a_{i_k j_k} -1 \right),$$
where the last sum is taken over all $(n!/(n-k)!)^2 \leq n^{2k}$ pairs of ordered sets $(i_1, \ldots, i_k)$ and $(j_1, \ldots, j_k)$ of distinct numbers between $1$ and $n$.
\hfill $\square$

It follows that the polynomial $p$ of Theorem \ref{th1.5} can be computed in time $n^{O(\ln n -\ln \epsilon)}$, where the implied constant in the ``$O$" notation  
depends on $\eta$ alone.

\subsection{Proof of Theorem \ref{th2.5}}\label{subsec7.4}
The proof is very similar to that of Section \ref{subsec7.3}. Let $J=J_{2n}$ be the $2n \times 2n$ matrix filled with 1s and let $A=\left(a_{ij}\right)$ be a $2n \times 2n$ symmetric complex matrix satisfying the conditions of theorem.
We define a univariate polynomial
$$g(z) =\haf \bigl(J + z(A-J)\bigr) \quad \text{for} \quad z \in {\mathbb C},$$
so that $\deg g \leq n$, 
$$g(0) =\haf J =\frac{(2n!)}{2^n n!} \quad \text{and} \quad g(1)=\haf A.$$
Moreover, by Theorem \ref{th2.4}, we have
$$g(z) \ne 0 \quad \textrm{provided} \quad |z| \leq \beta \quad \text{where} \quad \beta=\frac{0.5}{\eta} >1.$$
 We write
$$g(z)=\sum_{\{i_1,j_1\}, \ldots, \{i_n, j_n\}} \left(1 +z \left(a_{i_1j_1}-1\right)\right) \cdots \left(1+ z\left(a_{i_n j_n}-1\right) \right),$$
where the sum is taken over all $(2n)!/2^n n!$ unordered partitions of the set $\{1, 2, \ldots, 2n\}$ into $n$ pairwise disjoint unordered pairs 
$\{i_1, j_1\}, \ldots, \{i_n, j_n\}$. Hence 
for $k > 0$ we have
$$g^{(k)}(0)=\frac{k! (2n-2k)! }{2^{(n-k)} (n-k)!} \sum_{\{i_1, j_1\}, \ldots, \{i_k, j_k\}} \left(a_{i_1 j_1}-1\right) \cdots \left(a_{i_k j_k} -1 \right),$$
where the sum is taken over all unordered collections $\{i_1, j_1\}, \ldots, \{i_k, j_k\}$ of pairwise disjoint unordered pairs.

The proof then proceeds as in Section \ref{subsec7.3}.
\hfill $\square$

It follows that the polynomial $p$ of Theorem \ref{th2.5} can be computed in time $n^{O(\ln n-\ln \epsilon)}$, where the implied constant in the 
``$O$" notation depends on $\eta$ alone.

\subsection{Proof of Theorem \ref{th3.3}}\label{subsec7.5}
 Let $J=J_{n,d}$ be the $d$-dimensional $n \times \ldots \times n$ tensor filled with 1s and let $A$ be the tensor satisfying the conditions of the theorem. We introduce a univariate polynomial 
 \begin{equation}\label{eq7.5.1}
g(z)=\PER \bigl(J+ z(A-J)\bigr),
\end{equation}
so that 
$$g(0)=\PER J = (n!)^{d-1} \quad \text{and} \quad g(1)=\PER A.$$
Moreover, by Theorem \ref{th3.2},
$$g(z) \ne 0 \quad \text{provided} \quad |z| \leq \beta \quad \text{where} \quad \beta =\frac{\eta_d}{\eta} > 1.$$
We write
$$g(z) = \sum_{\sigma_2, \ldots, \sigma_d \in S_n} \prod_{i=1}^n \left(1 + z \left(a_{i \sigma_2(i) \ldots \sigma_d(i)} -1 \right)\right),$$
so that 
$$g(0)=\PER J = (n!)^{d-1}$$
and for $k >0$,
$$g^{(k)}(0)=\sum_{\sigma_2, \ldots, \sigma_d \in S_n} \sum_{(i_1, \ldots, i_k)}\left(a_{i_1 \sigma_2(i_1) \ldots \sigma_d(i_1)} -1 \right) 
\cdots \left(a_{i_k \sigma_2(i_k) \ldots \sigma_d(i_k)} -1 \right), $$
where the last sum is taken over all ordered $k$-tuples $(i_1, \ldots, i_k)$ of distinct indices $1 \leq i_j \leq n$.
By symmetry we can write
\begin{equation*}
\begin{split}
g^{(k)}(0)=&\left((n-k)!\right)^{d-1} \\ &\quad \times \sum_{\substack{ \left(i_{11}, \ldots i_{k1}\right) \\ \left(i_{12}, \ldots, i_{k2}\right) \\ \ldots \ldots \ldots \ldots \ldots \\
\left(i_{d1}, \ldots, i_{kd}\right)}} \left(a_{i_{11} i_{12} \ldots i_{1d}}-1 \right)\left(a_{i_{21} i_{22} \ldots i_{2d}}-1 \right)  \cdots \left(a_{i_{k1} i_{k2} \ldots i_{kd}}-1\right), \end{split}
\end{equation*}
where the last sum is taken over all $(n!/(n-k)!)^d \leq n^{kd}$ collections of $d$ ordered $k$-tuples $\left(i_{1j}, \ldots, i_{kj}\right)$ for $j=1, \ldots, d$ of distinct 
indices $1 \leq i_{1j}, \ldots, i_{kj} \leq n$. The proof then proceeds as in Section \ref{subsec7.3}.
\hfill $\square$

The polynomial $p$ can be computed in $n^{O(\ln n - \ln \epsilon)}$ time, where the implied constant in the ``$O$" notation depends on $\eta$ and 
$d$ only.

\subsection{Proof of Theorems \ref{th1.7} and \ref{th3.7}}\label{subsec7.6}
 As Theorem \ref{th1.7} is a particular case of Theorem \ref{th3.7}, we prove the latter theorem only. 
As in Section \ref{subsec7.5}, we define the univariate polynomial (\ref{eq7.5.1}). By Theorem \ref{th3.6}, we have 
$$g(z) \ne 0 \quad \text{provided} \quad |z| \ \leq \ \beta \quad \text{where} \quad \beta = \frac{\eta_d}{\eta}  \ > \ 1,$$
and the proof follows as in Section \ref{subsec7.5}. 
\hfill $\square$

\section{Proofs of Theorems \ref{th1.2}, \ref{th2.2} and \ref{th3.4}}\label{sec8}

Lemma \ref{le7.1} allows us to approximate the value of $\ln g(1)$ by a low degree Taylor polynomial of $\ln g(z)$ at $z=0$ provided the polynomial $g(z)$ does not have zeros in a disc of radius $\beta >1$ centered at $z=0$. In view of Theorems \ref{th1.3}, \ref{th2.3} and \ref{th3.5}, we would like to construct a similar approximation under a weaker assumption that $g(z) \ne 0$ for $z$ in some neighborhood of the interval $[0, 1]$ in the complex plane.
 To achieve that, we first construct a polynomial $\phi$ such that $\phi(0)=0$, $\phi(1)=1$ and such that $\phi$ maps the disc $|z| \leq \beta$ for some $\beta >1$ inside the neighborhood. We then apply Lemma \ref{le7.1} to the composition $g(\phi(z))$. The following lemma provides an explicit construction of such a polynomial $\phi$.

\begin{lem}\label{le8.1} For $0 < \rho < 1$, let us define
\begin{equation*}
\begin{split} \alpha=&\alpha(\rho)=1-e^{-\frac{1}{\rho}}, \quad \beta=\beta(\rho)=\frac{1 - e^{-1-\frac{1}{\rho}}}{ 1- e^{ -\frac{1}{\rho}}} \ > \ 1,\\
N=&N(\rho) = \left\lfloor \left(1 +\frac{1}{\rho}\right) e^{1 + \frac{1}{\rho}} \right\rfloor \ \geq \ 14, \quad 
\sigma=\sigma(\rho)=\sum_{m=1}^N \frac{\alpha^m}{m} \quad \text{and} \\
\phi(z)=&\phi_{\rho}(z)=\frac{1}{\sigma} \sum_{m=1}^N \frac{(\alpha z)^m}{m}.
\end{split}
\end{equation*}
Then $\phi(z)$ is a polynomial of degree $N$ such that $\phi(0)=0$, $\phi(1)=1$, 
$$-\rho \ \leq \  \Re\thinspace \phi(z) \ \leq \ 1+2 \rho \quad \text{and} \quad \left| \Im\thinspace \phi(z) \right| \ \leq \ 2 \rho 
\quad \text{provided} \quad |z| \leq \beta.$$
\end{lem}
\textbf{Proof.} Clearly, $\phi(z)$ is a polynomial of degree $N$ such that $\phi(0)=0$ and $\phi(1)=1$. It remains to prove that $\phi$ maps the 
disc $|z| \leq \beta$ into the strip $ -\rho \leq \Re\thinspace z \leq 1+2\rho$, $\left| \Im\thinspace z \right| \leq 2\rho$.

We consider the function
$$F_{\rho}(z)=\rho \ln \frac{1}{1-z} \quad \text{for} \quad |z| < 1.$$
Since 
$$\Re\thinspace \frac{1}{1-z} \ > \ \frac{1}{2} \quad \text{if} \quad |z| < 1, $$
the function $F_{\rho}(z)$ is well-defined by the choice of a branch of the logarithm, which we choose so that 
$$F_{\rho}(0)=\rho \ln 1 =0.$$
Then for $|z| < 1$ we have
\begin{equation}\label{eq8.1.1}
\left| \Im\thinspace F_{\rho}(z) \right| \ \leq \ \frac{\pi \rho}{2} \quad \text{and} \quad
\Re\thinspace F_{\rho}(z) \ \geq \ -\rho \ln 2. 
\end{equation}
In addition,
\begin{equation}\label{eq8.1.2}
F_{\rho}(\alpha)=1 \quad \text{and} \quad \Re\thinspace F_{\rho}(z) \ \leq \ 1+\rho  \quad \text{provided} \quad |z| \leq 1 -e^{-1-\frac{1}{\rho}}.
\end{equation}
Let 
$$P_n(z)=\sum_{m=1}^n \frac{z^m}{m}.$$
Then
$$\left| \ln \frac{1}{1-z} -P_n(z)\right| =\left| \sum_{m=n+1}^{\infty} \frac{z^m}{m} \right| \ \leq \ \frac{|z|^{n+1}}{(n+1)(1-|z|)} \quad \text{provided} \quad 
|z| < 1.$$
Therefore, for $|z| \leq \beta$, we have
\begin{equation}\label{eq8.1.3}
\begin{split}
\left| F_{\rho} (\alpha z) - \rho P_N(\alpha z) \right| \ \leq \ &{\rho} \frac{(\alpha \beta)^{N+1}}{(N+1)(1-\alpha \beta)} \\ = \ 
&\frac{\rho}{N+1} \left(1 - e^{-1 -\frac{1}{\rho}}\right)^{N+1} e^{1+ \frac{1}{\rho}} \\ \leq \ &\frac{\rho}{N+1} \ \leq \ \frac{\rho}{15}.
\end{split}
\end{equation}
Combining (\ref{eq8.1.1}), (\ref{eq8.1.2}) and (\ref{eq8.1.3}), we conclude that for $|z| \leq \beta$ we have 
\begin{equation}\label{eq8.1.4}
\left| \Im\thinspace \rho P_N(\alpha z) \right| \ \leq \ 1.64 \rho \quad \text{and} \quad -0.76 \rho \ \leq \ \Re\thinspace \rho P_N(\alpha z) \ \leq \ 1+ 1.07\rho. 
\end{equation}
Substituting $z=1$ in (\ref{eq8.1.3}) and using (\ref{eq8.1.2}), we conclude that 
\begin{equation}\label{eq8.1.5}
\left|1 - \rho P_N(\alpha) \right| \ \leq \ \frac{\rho}{15}. 
\end{equation}
We have
$$\phi(z)= \frac{P_N(\alpha z)}{P_N(\alpha)} = \frac{\rho P_N(\alpha z)}{\rho P_N(\alpha)},$$
where $\rho P_N(\alpha)$ is positive real, which by (\ref{eq8.1.5}) satisfies 
$$\rho P_N(\alpha) \ \geq \ 1 -\frac{\rho}{15}.$$
Since
$$\left(1-\frac{\rho}{15}\right)^{-1} \ \leq \ \min\left\{\frac{15}{14},\ 1+ \frac{2 \rho}{15} \right\} \quad \text{for} \quad 0 \leq \rho \leq 1,$$
from (\ref{eq8.1.4}) we conclude that 
$$\left| \Im\thinspace \phi(z) \right| \ \leq \ 2 \rho \quad \text{and} \quad -\rho \ \leq \ \Re\thinspace \phi(z) \ \leq \  1+ 2\rho \quad 
\text{provided} \quad |z| \leq \beta.$$
\hfill $\square$

\subsection{Proof of Theorem \ref{th1.2}}\label{subsec8.2} Let $A=\left(a_{ij}\right)$ be an $n \times n$ real matrix satisfying the conditions of the theorem and
let $J=J_n$ be the $n \times n$ matrix filled with 1s. As in Section \ref{subsec7.3}, we define a univariate polynomial 
$$r(z)=\per\bigl(J +z(A-J)\bigr) \quad \text{for} \quad z \in {\mathbb C}. $$
Suppose that 
\begin{equation}\label{eq8.2.1}
-\xi \ \leq \ \Re\thinspace z \ \leq \ 1+\xi \quad \textrm{and} \quad \left| \Im\thinspace z \right| \ \leq \ \zeta 
\end{equation}
for some $\xi>0$ and $\zeta >0$. Then the entries $b_{ij}$ of the matrix $B=J+z(A-J)$ satisfy 
$$\left| 1-\Re\thinspace b_{ij}\right| \ \leq \ (1+\xi)(1-\delta) \quad \text{and} \quad \left| \Im\thinspace b_{ij} \right| \ \leq \ \zeta(1-\delta).$$
We choose $\xi=\xi(\delta) >0$ such that 
$$\eta=(1+\xi)(1-\delta) \ < \ 1$$ and then choose 
$\zeta=\zeta(\delta) > 0$ such that 
$$\zeta(1-\delta) \ < \ (1-\eta) \sin \left( \frac{\pi}{4} - \arctan \eta \right).$$
By Theorem \ref{th1.3}, we have $r(z) \ne 0$ for $z$ satisfying (\ref{eq8.2.1}). 

Using Lemma \ref{le8.1}, we construct a univariate polynomial $\phi$ of some degree $N=N(\delta)$ such that $\phi(0)=0$, $\phi(1)=1$ and 
$\phi$ maps the disc $\{z: \ |z| \leq \beta \}$ inside the strip (\ref{eq8.2.1}), where $\beta=\beta(\delta) > 1$.
Let  
$$g(z)=r(\phi(z)) \quad \textrm{for} \quad z \in {\mathbb C}.$$
Then $g(z)$ is a univariate polynomial such that $\deg g \ \leq \ N n$, 
$$g(0)=r(0)=\per J =n! \quad \textrm{and} \quad g(1)=r(1)=\per A.$$
Besides,
$$g(z) \ne 0 \quad \textrm{provided} \quad |z| \leq \beta.$$ Let us define 
$$f(z) =\ln g(z) \quad \text{for} \quad |z| \leq \beta,$$
where we chose the branch of the logarithm such that $f(0)=\ln n!$ is real. Let 
$T_m(z)$ be the Taylor polynomial of $f(z)$ of degree $m$ computed at $z=0$. By Lemma \ref{le7.1}, we have 
$$\left|T_m(1)-f(1)\right| =\left| T_m(1) - \ln \per A \right| \ \leq \ \epsilon,$$
for some $m \leq \gamma(\ln n - \ln \epsilon)$ where $\gamma=\gamma(\delta)>0$ is a constant depending on $\delta$ alone.
It remains to show that $T_m(1)$ is a polynomial in the entries of $A$ of degree not exceeding $m$.

For a univariate polynomial $p(z)$, let $p_{[m]}$ be the polynomial obtained from $p$ by discarding all monomials of degree higher than $m$.
 Since $\phi(0)=0$, the constant term of of $\phi$ is $0$ and 
therefore 
$$g_{[m]}=\left(r(\phi)\right)_{[m]} = \left(r_{[m]}\left(\phi_{[m]}\right)\right)_{[m]}.$$
In words: to compute the polynomial $g_{[m]}$ obtained from $g$ by discarding the monomials of degree higher than $m$, it suffices to compute 
the polynomials $r_{[m]}$ and $\phi_{[m]}$ obtained from $r$ and $\phi$ respectively by discarding the monomials of degree higher than $m$, and then discard the monomials of degree higher than $m$ in the composition $r_{[m]}(\phi_{[m]})$.

From Section \ref{subsec7.3}, it follows that $r^{(k)}(0)$ is a polynomial of degree $k$ in the entries of the matrix $A$. It follows then that $g^{(k)}(0)$ is a polynomial in the entries of $A$ of degree at most $k$ that can be computed in $n^{O(k)}$ time (the implied constant in the ``$O$'' notation is absolute).
From Section \ref{subsec7.2} it follows then that $f^{(k)}(0)$ is a polynomial in the entries of $A$ of degree at most $k$, which completes the proof.
\hfill $\square$

\subsection{Proof of Theorem \ref{th2.2}}\label{subsec8.3} Given a $2n \times 2n$ real symmetric matrix $A$ satisfying the conditions of the theorem, 
we define the univariate polynomial $r(z)$ by 
$$r(z)=\haf\bigl(J+z(A-J)\bigr),$$
where $J=J_{2n}$ is the $2n \times 2n$ matrix filled with 1s and the proof then proceeds as in Section \ref{subsec8.2}, only that the reference to Theorem \ref{th1.3} is replaced by the reference to Theorem \ref{th2.3} and the reference to Section \ref{subsec7.3} is replaced by the reference to Section \ref{subsec7.4}.
\hfill $\square$

\subsection{Proof of Theorem \ref{th3.4}}\label{subsec8.4}Given a $d$-dimensional $n \times \ldots \times n$ tensor $A$ satisfying the conditions of the theorem, 
we define the univariate polynomial $r(z)$ by 
$$r(z) = \PER (J+ z(A-J)),$$
where $J=J_{d,n}$ is the $d$-dimensional $n \times \ldots \times n$ tensor filled with 1s. Suppose that (\ref{eq8.2.1}) holds for some $\xi >0$ and $\zeta > 0$. 
Then the entries $b_{i_1 \ldots i_d}$ of the tensor $B=J+z(A-J)$
satisfy 
$$\left|1 - \Re\thinspace b_{i_1 \ldots i_d} \right| \ \leq \ (1+\xi) \eta \quad \text{and} \quad \left| \Im\thinspace b_{i_1 \ldots i_d} \right| \ \leq \ \zeta \eta.$$
We choose $\xi=\xi(\eta) >0$ such that 
$$\eta'=(1+\xi) \eta \ < \ \tan \frac{\pi}{4(d-1)}$$ 
and then choose $\zeta=\zeta(\eta) > 0$ such that 
$$\zeta \eta \ < \ (1-\eta') \sin\left(\frac{\pi}{4(d-1)}  -\arctan \eta' \right).$$
By Theorem \ref{th3.5}, we have $r(z) \ne 0$ for $z$ satisfying (\ref{eq8.2.1}) with $\xi$ and $\zeta$ so chosen. The proof then proceeds as in Section \ref{subsec8.2}, only that the reference to Section \ref{subsec7.3} is replaced by the reference to Section \ref{subsec7.5}.
\hfill $\square$

\section{Concluding remarks}\label{sec9}

\subsection{Numerical experiments}\label{subsec9.1}
The algorithm of Theorem \ref{th1.5} for approximating permanents of real and complex matrices was implemented by Kontorovich and Wu \cite{KW16}, who
conducted numerical experiments on approximating $\ln \per A$ by a polynomial of just degree 3. The experiments seem to show that the method a) very fast, b) quite accurate on positive matrices with entries within a factor of 10 of each other and c) quite accurate on random 0-1 matrices with at most 10$\%$ of zeros.

\subsection{Connections to the Szeg\H{o} curve}\label{subsec9.2}
Let $I=I_n$ be the $n \times n$ identity matrix, let $J=J_n$ be the $n \times n$ matrix of 1s and let 
$r(z)=\per (J +z(nI -J))$
 be the polynomial of Section \ref{subsec8.2} for the matrix $A=nI$. We have 
 \begin{equation*}
 \begin{split} r(z)=&\per \bigl( z n I + (1-z) J  \bigr)= (1-z)^n \per \left( \frac{z}{1-z} nI + J\right)\\=
 &(1-z)^n \sum_{k=0}^n \binom{n}{k} \left(\frac{z}{1-z}\right)^k n^k (n-k)! = (1-z)^n n! \sum_{k=0}^n \frac{(ny)^k}{k!}, \\
 &\qquad \textrm{where} \quad
 y=\frac{z}{1-z}. 
 \end{split}
 \end{equation*}
 Kontorovich and Wu noticed \cite{KW16} that the location of the complex zeros of $r(z)$, which is crucial for our analysis of the approximation of the permanent,
 for $A=nI$ can be determined from a result of Szeg\H{o}, who showed in 1922 that as $n \longrightarrow \infty$, the zeros of the polynomial
 $$\sum_{k=0}^n \frac{(ny)^k}{k!}$$ 
 converge to the curve $\left\{ \left| \zeta e^{1-\zeta}\right| =1,\ |\zeta| \leq 1 \right\}$, now known as the {\it Szeg\H{o} curve}, cf. \cite{Wa03}. It follows then that the roots of $r(z)$ in the vicinity of the interval $[0, 1]$ for large $n$ cluster around $z=0.5$, so our method of interpolation from $J$ to $A=nI$ works roughly ``halfway". The same is true if $A=nP$, where $P$ is an $n \times n$ permutation matrix, and there is some limited computational evidence 
 that for non-negative $n \times n$ matrices $A$ with row and column sums $n$ (such matrices are convex combinations of matrices $nP$) the polynomial 
 $r(z)=\per (J +z(A-J))$ has no zeros in the vicinity of the interval $[0, 0.5-\epsilon]$ for any $\epsilon >0$ and all sufficiently large $n$, so our method works ``at 
 least halfway" for all such matrices $A$, cf. also \cite{Mc14}. Combined with the scaling algorithm, see \cite{L+00}, this may lead to a useful algorithm for approximating permanents of arbitrary non-negative matrices.
 
 \subsection{Connections to the mixed characteristic polynomial}\label{subsec9.3}
 In their solution of the Kadison - Singer problem, Marcus, Spielman and Srivastava
 \cite{M+15} introduced and studied the {\it mixed characteristic polynomial} of $n$ Hermitian $n \times n$ matrices $Q_1, \ldots, Q_n$, 
 $$p_{Q_1, \ldots, Q_n}(x) =\prod_{i=1}^n \left(1 - \frac{\partial}{\partial z_i}\right) \det\left(xI + \sum_{i=1}^n z_i Q_i \right) \Big|_{z_1=\ldots = z_n=0},$$
 where $I$ is the $n \times n$ identity matrix. In particular, they showed that the roots of $p$ are necessarily non-negative real provided $Q_1, \ldots, Q_n$ are 
 non-negative semidefinite. An anonymous referee pointed out to a similarity between the mixed characteristic polynomial and the polynomial 
 $r(z)=\per (J +z (A-J))$ used in this paper. Given an $n \times n$ non-negative matrix $A$, let $Q_i$ be the diagonal matrix having the $i$-th row of $A$ as the diagonal. Then $Q_1, \ldots, Q_n$ are non-negative semidefinite matrices and the mixed characteristic polynomial can be written 
 as 
 $$p_A(-x)=(-1)^n \sum_{k=0}^n W_{n-k}(A) x^k,$$
 where $W_k(A)$ is the sum of permanents of the $k \times k$ submatrices of $A$, so up to a sign and a substitution $x \longmapsto -1/x$, the polynomial $p_A$ is the matching polynomial of Section \ref{sec6} (and the fact that the roots of $p_A$ are non-negative real is a particular case of the Heilmann - Lieb Theorem \cite{HL72}).
 On the other hand,
 $$r(z) = z^n \sum_{k=0}^n k! W_{n-k}(A) \left(\frac{1-z}{z}\right)^k.$$
 The relation between $p$ and $r$ is essentially used in the proof of Theorem \ref{th3.6}, which was absent in the version of the paper the referee commented on, 
 but was obtained before the author received the comment.
 
 On the other hand, the general mixed characteristic polynomial may appear useful for approximating the {\it mixed discriminant} of $Q_1, \ldots, Q_n$, which, up to a sign is just the constant term of $p_{Q_1, \ldots, Q_n}$.

 \subsection{Approximation of general polynomials}\label{subsec9.4}
  Lemmas \ref{le7.1} and \ref{le8.1} suggest the following general way of approximating combinatorially 
interesting polynomials. Suppose that $p(z)$ is a univariate polynomial such that $\deg p \leq n$. Suppose further we want to approximate 
$p(1)$ whereas $p(0)$ is easily computable and the derivatives $p^{(k)}(0)$ can be computed in $n^{O(k)}$ time. We can approximate $p(1)$ within 
a relative error $\epsilon>0$ in quasi-polynomial time $n^{O(\ln n -\ln \epsilon)}$ provided we can find a ``sleeve" $S \subset {\mathbb C}$ in the complex plane such that $0 \in S$, $1 \in S$ and $p(z) \ne 0$ for all $z \in S$. The sleeve $S$ should be wide enough, meaning that it contains a 
number $N$, fixed in advance, of discs $D_1, \ldots, D_N$ of equal radii such that $D_i$ contains the center of $D_{i-1}$ 
for $i=2, \ldots, N$ with $D_1$ centered at $0$ and $D_N$ centered at $1$. An example of such a sleeve is provided by the strip $-\delta \leq \Re\thinspace z \leq 1+\delta$ and $\left|\Im\thinspace z \right| \leq \tau$ for some $\delta >0$ and $\tau >0$, fixed in advance for the polynomial $r(z)=\per \left(J+z(A-J)\right)$ of 
Section \ref{subsec8.2}.

As another example, we consider the independence polynomial of graph. Let $G=(V, E)$ be a graph (undirected, without loops or multiple edges) with set $V$ of vertices and set $E$ of edges. A set $S \subset V$ is called {\it independent} if 
no two vertices in $S$ span an edge of $G$ (the empty set $S=\emptyset$ is considered independent). The {\it independence polynomial} of $G$ is 
defined as 
\begin{equation*}
p_G(z)=\sum_{\substack{S \subset V \\ S \text{\ is independent}}} z^{|S|} = \sum_{k=0}^{|V|} \left(\text{the number of independent $k$-sets in $V$}\right)z^k.
\end{equation*}
Then $p_G(1)$ is the number of all independent sets in $G$, a quantity of considerable combinatorial interest. On the other hand, 
the value of the derivative $p_G^{(k)}(0)$ can be computed in $|V|^{O(k)}$ time by a direct inspection of all $k$-subsets $S \subset V$.

Suppose we know that $p_G(z) \ne 0$ provided $|z| \leq \beta$ for some $\beta >0$ (for example, $\beta$ can be the Dobrushin bound, see \cite{SS05}
and \cite{CF16}). Lemma \ref{le7.1} then implies that for any $0 \leq \lambda < 1$, fixed in advance, the value of $p_G(z)$ can be approximated within a relative error $0 < \epsilon < 1$ in quasi-polynomial time $|V|^{O(\ln |V| - \ln \epsilon)}$ provided $|z| \leq \lambda \beta$, see \cite{Re15} for many 
examples of this nature and also \cite{We06} and \cite{H+16} for algorithms based on the ``correlation decay" idea.

If, additionally,  the zeros of $p_G(z)$ are known to be confined to a particular region of the complex plane ${\mathbb C}$, we can hope to do better by constructing 
a sleeve $S \subset {\mathbb C}$ where $p_G(z)$ is not zero and interpolating $p_G(z)$ there. In an extreme case, when $G$ is claw-free, the roots of 
$p_G(z)$  are known to be negative real \cite{CS07}, which leads to a quasi-polynomial algorithm for approximating $p_G(z)$ provided 
$\left|\pi -\arg\thinspace z\right| > \lambda^{-1}$ (so that $z$ stays away from the negative real axis) and 
$|z| \leq \lambda \beta$ where $\lambda >0$ is arbitrarily large, fixed in advance, see also \cite{B+07} for an algorithm based on the correlation decay approach.

On the other hand, for a general graph $G$ there cannot be such a sleeve $S$ unless NP-complete problems admit a quasi-polynomial time algorithm.
Indeed, generally, it is an NP-hard problem to approximate $p_G(z)$ for a real $z > \lambda \beta$, where $\lambda >0$ is some absolute constant and 
$\beta$ is the Dobrushin lower bound on the absolute value of the roots of $p_G(z)$ \cite{LV99}. This means that for a general graph $G$ one can expect the complex roots of $p_G(z)$ to ``surround" the origin, so that there is no possibility to squeeze a sleeve between them to connect 0 and 1. 

Since the first version of this paper appeared as a preprint, this general direction was pursued further in \cite{PR16}.

\subsection{Approximating multi-dimensional permanents better}\label{subsec9.5}
It would be interesting to extend the class of polynomials for which a version of Theorems \ref{th1.2} and \ref{th2.2} can be obtained. While we failed to obtain such a version for the multi-dimensional permanent (see Section \ref{sec3}), there does not seem to be a computational complexity obstacle for such an extension to exist.
In \cite{BS11} it is shown that the $d$-dimensional permanent of a $n \times \ldots \times n$ tensor with positive entries between an arbitrarily small $\delta >0$, fixed in advance, and 1 can be approximated within an 
$n^{O(1)}$ factor in polynomial time, where the implicit constant in the ``$O$" notation depends only on $d$ and $\delta$, which can be viewed as an indirect evidence that Theorem \ref{th1.2} can indeed be extended to multi-dimensional permanents.

\section*{Acknowledgments}

I am grateful to anonymous referees for their careful reading of the paper and suggestions and to Max Kontorovich and Han Wu for conducting numerical experiments on the approximation of permanents and pointing out to connections with the Szeg\H{o} curve.

\bibliographystyle{amsplain}

\begin{dajauthors}
\begin{authorinfo}[abar]
  Alexander Barvinok\\
  Professor\\
  University of Michigan \\
  Ann Arbor, MI, USA\\
  barvinok\imageat{}umich\imagedot{}edu \\
  \url{http://www.math.lsa.umich.edu/~barvinok/}
\end{authorinfo}
\end{dajauthors}

\end{document}